# An improved characterisation of regular generalised functions of white noise and an application to singular SPDEs

**Martin Grothaus[1]** · **Jan Müller[1]** · **Andreas Nonnenmacher[1]**

## Abstract

A characterisation of the spaces $\mathcal{G}_K$ and $\mathcal{G}'_K$ introduced in Grothaus et al. (Methods Funct Anal Topol 3(2):46–64, 1997) and Potthoff and Timpel (Potential Anal 4(6):637–654, 1995) is given. A first characterisation of these spaces provided in Grothaus et al. (Methods Funct Anal Topol 3(2):46–64, 1997) uses the concepts of holomorphy on infinite dimensional spaces. We, instead, give a characterisation in terms of U-functionals, i.e., classic holomorphic function on the one dimensional field of complex numbers. We apply our new characterisation to derive new results concerning a stochastic transport equation and the stochastic heat equation with multiplicative noise.

The third author thanks the department of Mathematics at the University of Kaiserslautern for financial support in the form of a fellowship.

✉ Martin Grothaus
  grothaus@mathematik.uni-kl.de

  Jan Müller
  jan.mueller.mathematik@gmail.com

  Andreas Nonnenmacher
  nonnenmacher@mathematik.uni-kl.de

[1] Mathematics Department, TU Kaiserslautern, PO Box 3049, 67653 Kaiserslautern, Germany

�underline 2 Springer

## 1 Introduction

Gaussian Analysis, in particular White Noise Analysis, has been intensively investigated and developed in recent years. It gained interest by its applications in stochastic (partial) differential equations, quantum physics and many more. One aspect is the construction, analysis and characterisation of spaces of (generalised) random variables on infinite dimensional Gaussian spaces. In this paper we deal with a specific type of random variables, which are natural in the context of stochastic partial differential equations, as illustrated in Sect. 5 below via different examples. These random variables have important properties, such as Malliavin differentiability. The type of random and generalised random variables under consideration in this paper can be described as follows. Let $\mathcal{H}$ be a real seperable Hilbert space carrying a self-adjoint operator $(K, D(K))$. Furthermore, $\mathcal{N}$ is a nuclear space, densely and continuously embedded into $\mathcal{H}$ such that $\mathcal{N} \subset D(K)$. By the Bochner-Minlos theorem one obtains a Gaussian measure $\mu$ on the dual space $\mathcal{N}'$ of $\mathcal{N}$ with covariance functional $(\cdot, \cdot)_{\mathcal{H}}$. Via the Wiener–Itô–Segal isomorphism the second quantisation $\Gamma(K)$ of $(K, D(K))$ is defined on the space $L^2(\mathcal{N}', \mu)$. The random variables $\mathcal{G}_K \subset L^2(\mathcal{N}', \mu)$ we investigate are exactly the $C^\infty$ vectors of the self-adjoint operator $\Gamma(K)$. Furthermore, the operator $\Gamma(K)$ induces a finer topology on $\mathcal{G}_K$. The space of generalised random variables $\mathcal{G}_K'$ is the dual space w.r.t. this topology. Important examples of random variables and their dual space arise in this way. For example, the pair of Hida test functions and distributions $(S)$ and $(S)'$, see e.g. [15], and the pair $\mathcal{G}$ and $\mathcal{G}'$ in [25] arise in this way for suitable choices of the operator $K$. In particular, for $K = \lambda Id, \lambda > 1$, the elements of the space $D(\Gamma(K))$, which contains $\mathcal{G}_K$, are infinitely often Mallivain differentiable along $H$.

Our results can be divided into two parts. The first part consists of a refinement of the one found in [12]. The authors used the concept of holomorphy on Hilbert spaces. In this paper we avoid this technique, which also results in a shorter proof of the main result. Furthermore, this makes our result easier to apply. To overcome the usage of holomorphy on Hilbert spaces we use the concepts of the $S$-transform (see Definition 2.6) and $U$-functionals (see Definition 2.7) as well as the famous characterisation theorem by Potthoff and Streit (see Theorem 2.8). Observe that in applications (generalised) random variables are often constructed and defined only via their $S$-transform, see also the example in Sect. 5 below. Fortunately, this is the only ingredient we need for our characterisation. In the second part we deal with two different kind of stochastic partial differential equations. The first one is a stochastic transport equation, the second is the stochastic heat equation both with a multiplicative noise. For both equations we give explicit conditions in terms of the coefficients of the equations such that their respective solutions are actually contained in the much smaller space $\mathcal{G}_{aI} \subseteq L^2(\mu), a > 1$, of smooth functions in the sense of Malliavin calculus.

This article is organized as follows. In Sect. 2 we briefly describe the functional analytic framework and the main concept of Gaussian and white noise analysis we use throughout this paper to state our main theoretical result. In particular we give the definition of the spaces $\mathcal{G}_K$ and $\mathcal{G}_K'$ under consideration. Theorem 2.11 contains our main result, a characterisation of $\mathcal{G}_K$ and $\mathcal{G}_K'$ in terms of the $S$-transform. In Sect.

3 we further introduce concepts of Gaussian Analysis. Section 4 contains the proof of Theorem 2.11. In Sect. 5 we present two applications from the field of stochastic partial differential equations. In the first case we apply our main result to the stochastic partial differential equation

$$\left.\begin{array}{l} \dfrac{\partial u_{t,x}}{\partial t} = \dfrac{1}{2}\nu(t)\dfrac{\partial^2 u_{t,x}}{\partial x^2} + \dfrac{\partial u_{t,x}}{\partial x}\sigma(t)\dot{B}_t \quad t > 0, x \in \mathbb{R}, \\ u(0, \cdot) = \delta_0, \end{array}\right\} (STE)$$

where $\frac{\partial u_{t,x}}{\partial x}\sigma(t)\dot{B}_t$ is understood in the Itô sense. The coefficients $\nu$ and $\sigma$ are allowed to be singular. This equation was treated by several authors, see e.g. [4,11,24]. In particular, in [24] solutions were constructed as elements of the Hida distribution space $(S)'$, see the references for the precise statement. We use the characterisation theorem to improve the results from [24] by showing that the solution belongs to the space of regular distributions $\mathcal{G}'$. In particular, we determine explicitly in terms of the coefficients the regularity of the solutions, see Theorem 5.2.

In the second part we consider a stochastic heat equation with general coloured noise, i.e.,

$$\left.\begin{array}{l} \dfrac{\partial u_{t,x}}{\partial t} = \dfrac{1}{2}\Delta u_{t,x} + u_{t,x}\dot{W}_{t,x}, \quad t > 0, x \in \mathbb{R}^d, \\ u_{0,x} = u_0(x), \quad x \in \mathbb{R}^d, \end{array}\right\} (SHE)$$

where the product between $u_{t,x}$ and the centered Gaussian process $\dot{W}_{t,x}$, $t > 0$, $x \in \mathbb{R}^d$, is treated in the Skorokhod and Stratonovich sense. The covariance of $\dot{W}_{t,x}$ is given in (31) below. Our results are based on [17] and extend the results given there. In particular, we show that $u(t, x) \in \mathcal{G} \subseteq L^2(\mu)$ for all $t \in (0, \infty), x \in \mathbb{R}^d$. This implies, using the results from [25], that $u(t, x)$ is infinitely often Malliavin differentiable. This was not shown in [17]. Eventually, in Sect. 6 we give an outlook for further applications of the derived characterisation in the context of stochastic currents.

The following core results are achieved in this article:

(i) We prove a new characterisation theorem for the space $\mathcal{G}_K$ and its dual $\mathcal{G}'_K$, which is an improvement of the result in [12].

(ii) In Example 2.12 we show how to construct appropriate nuclear triples to use our theoretical result in (i) in order to analyse stochastic partial differential equations driven by a Gaussian noise.

(iii) We derive explicit integrability conditions on the coefficients $\nu$ and $\sigma$ of (STE) to determine that the solution $u_{t,x}$ belongs to $\mathcal{G}'_{\lambda Id}, \lambda > 0$.

(iv) For the Skorokhod and Stratonovich version of (SHE) we improve results obtained in [17] and show that the corresponding mild solution is contained in $\mathcal{G}_{\lambda Id}, \lambda > 0$. This implies that the solution is smooth in the sense of Malliavin calculus.

The aim of this article is to further bridge the gap between classical stochastic analysis and white noise analysis. Moreover, it is intended to show case that the combination of white noise analysis and Malliavin calculus can be very fruitful.

## 2 Preliminaries and main results

To state our results we briefly introduce the main concepts of Gaussian analysis. The material in the following can be found in e.g. [2,15,19,22]. Henceforth in the Sects. 3, 4 and 5 we fix a separable real Hilbert space $(\mathcal{H}, (\cdot, \cdot)_{\mathcal{H}})$. Furthermore, there exists a real nuclear countably Hilbert space $\mathcal{N}$ densely and continuously embedded into $\mathcal{H}$. In the following we briefly explain the notion of a nuclear countably Hilbert space. I.e., there exists a family of real inner products $\{(\cdot, \cdot)_p\}_{p \in \mathbb{N}_0}$ on $\mathcal{N}$ with induced norms $\{\|\cdot\|_p\}_{p \in \mathbb{N}_0}$, where $(\cdot, \cdot)_0 = (\cdot, \cdot)_{\mathcal{H}}$. Theses norms satisfy $\|\varphi\|_p \leq \|\varphi\|_{p+1}$ for all $\varphi \in \mathcal{N}$ and $p \in \mathbb{N}_0$. Furthermore the family $\{\|\cdot\|_p\}_{p \in \mathbb{N}_0}$ is compatible, meaning that for all $p, q \in \mathbb{N}_0$ and every sequence $(\varphi_n)_{n \in \mathbb{N}} \subseteq \mathcal{N}$ which is a fundamental sequence w.r.t. $\|\cdot\|_q$ and converges to zero w.r.t. $\|\cdot\|_p$ converges also to zero w.r.t. $\|\cdot\|_q$. This implies that the identity operator $I : (\mathcal{N}, \|\cdot\|_p) \longrightarrow (\mathcal{N}, \|\cdot\|_q)$, $p > q$ extends linearly to an continuous, injective map with dense range from $\mathcal{H}_p$ to $\mathcal{H}_q$, where $\mathcal{H}_p$ and $\mathcal{H}_q$ denote the completion of $\mathcal{N}$ w.r.t. $\|\cdot\|_p$ and $\|\cdot\|_q$, respectively. This extension is denoted by $I_{p,q}$. Also, for every $q \in \mathbb{N}$ there exists a $p \geq q$ s.t. $I_{p,q}$ is a Hilbert-Schmidt operator. Eventually, the space $\mathcal{N}$ equipped with the metric $d(\varphi, \psi) = \sum_{p=0}^{\infty} 2^{-p} \frac{\|\varphi - \psi\|_p}{1 + \|\varphi - \psi\|_p}$ is assumed to be a seperable complete metric space. Hence, we obtain a chain of continuous and dense embeddings

$$\mathcal{N} \subseteq \mathcal{H}_p \subseteq \mathcal{H}_q \subseteq \mathcal{H} \subseteq \mathcal{H}_{-p} \subseteq \mathcal{H}_{-q} \subseteq \mathcal{N}', \quad p \geq q \tag{1}$$

where $\mathcal{N}'$, $\mathcal{H}_{-p}$ and $\mathcal{H}_{-q}$ denote the dual spaces of $\mathcal{N}$, $\mathcal{H}_p$ and $\mathcal{H}_q$, respectively. The dual pairing between an element $\varphi \in \mathcal{N}$ and $\Phi \in \mathcal{N}'$ is denoted by $\langle \varphi, \Phi \rangle \in \mathbb{R}$. We consider $\mathcal{N}'$ to be equipped with the weak topology and denote the respective Borel $\sigma$-field by $\mathcal{F}$. Via the Bochner-Minlos theorem we obtain measures defined on $\mathcal{N}'$ in the following way:

**Definition 2.1** Let $\sigma^2 > 0$ and define the continuous function

$$C_{\sigma^2} : \mathcal{N} \longrightarrow \mathbb{C}, \varphi \mapsto \exp\left(-\frac{\sigma^2}{2}(\varphi, \varphi)_{\mathcal{H}}\right). \tag{2}$$

Observe that $C_{\sigma^2}$ is positive definite and satisfies $C_{\sigma^2}(0) = 1$. Hence, by the Bochner-Minlos theorem, see e.g. [22, Theorem 1.5.2], we obtain a probability measure $\mu_{\sigma^2}$ defined on the Borel $\sigma$-field $\mathcal{F}$ of $\mathcal{N}'$ uniquely determined by the characteristic function $C_{\sigma^2}$, i.e., it holds

$$\int_{\mathcal{N}'} \exp\left(i \langle \varphi, \cdot \rangle\right) d\mu_{\sigma^2} = C_{\sigma^2}(\varphi) \quad \text{for all } \varphi \in \mathcal{N}. \tag{3}$$

For $\sigma^2 = 1$ we simply write $\mu$ instead of $\mu_1$. We denote by $L^2(\mu) := L^2(\mathcal{N}', \mathbb{C}; \mu)$ the space of equivalence classes of complex-valued functions which are square-integrable with respect to $\mu$. The next proposition is an immediate consequence of (2) and (3).

**Proposition 2.2** *Let $\varphi_1, ..., \varphi_n \in \mathcal{N}$, $n \in \mathbb{N}$. The image measure of $\mu$ under the map*

$$T_{\varphi_1, ..., \varphi_n} : \mathcal{N}' \longrightarrow \mathbb{R}^n, \omega \mapsto (\langle \varphi_i, \omega \rangle)_{i=1,...,n}$$

*is the Gaussian measure with mean zero and covariance matrix $C = \left((\varphi_i, \varphi_j)_{\mathcal{H}}\right)_{1 \leq i, j \leq n}$ on $\mathbb{R}^n$, i.e.,*

$$\mu \circ T_{\varphi_1, ..., \varphi_n}^{-1} = N(0, C).$$

An important subspace of $L^2(\mu)$ is the space of polynomials $\mathcal{P}(\mathcal{N}')$ on $\mathcal{N}'$. A polynomial $F \in \mathcal{P}(\mathcal{N}')$ is a function on $\mathcal{N}'$ of the form $F(\omega) = p(\langle \varphi_1, \omega \rangle, ..., \langle \varphi_k, \omega \rangle))$, where $k \in \mathbb{N}$, $\omega \in \mathcal{N}'$, $\varphi_1, ..., \varphi_k \in \mathcal{N}$, and $p$ is a complex polynomial in $k$ variables. An elementary proof shows that $\mathcal{P}(\mathcal{N}')$ is dense in $L^2(\mu)$. The subspace $\mathcal{P}_{(n)}(\mathcal{N}')$, $n \in \mathbb{N}_0$, is the space of all polynomials $F$ where $p$ is of degree at most $n$. Now we define orthogonal subspaces $W_{(n)}(\mu)$ for $n \in \mathbb{N}_0$. Define

$$W_{(0)}(\mu) := \text{span}(1),$$
$$W_{(n)}(\mu) := \mathcal{P}_{(n-1)}(\mathcal{N}')^{\perp} \cap \overline{\mathcal{P}_{(n)}(\mathcal{N}')}, \quad n \in \mathbb{N},$$

where $\overline{\mathcal{P}_{(n)}(\mathcal{N}')}$ denotes the closure of $\mathcal{P}_{(n)}(\mathcal{N}')$ and $\mathcal{P}_{(n-1)}(\mathcal{N}')^{\perp}$ the orthogonal complement of $\mathcal{P}_{(n-1)}(\mathcal{N}')$ in $L^2(\mu)$, respectively. For $n \in \mathbb{N}$ the subspace $W_{(n)}$ is called the space of $n$th order chaos. It follows by definition of $\left(W_{(n)}(\mu)\right)_{n \in \mathbb{N}_0}$ and the density of $\mathcal{P}(\mathcal{N}')$ in $L^2(\mu)$ that $L^2(\mu)$ is the orthogonal sum of the subspaces $W_{(n)}(\mu)$, $n \in \mathbb{N}_0$. A characteristic element of $W_{(n)}(\mu)$, $n \in \mathbb{N}_0$, is given by

$$\mathcal{N}' \ni \omega \mapsto H_{n, (\varphi, \varphi)_{\mathcal{H}}} (\langle \varphi, \omega \rangle) \in \mathbb{R},$$

where $\varphi \in \mathcal{N}$ and $H_{n, (\varphi, \varphi)_{\mathcal{H}}}$ is the $n$th Hermite polynomial with parameter $(\varphi, \varphi)_{\mathcal{H}}$. The family of Hermite polynomials with parameter $\alpha^2 > 0$ is defined via its generating function

$$\exp\left(-\alpha^2 \frac{t^2}{2} + tx\right) = \sum_{n=0}^{\infty} H_{n, \alpha^2}(x) \frac{t^n}{n!}, \quad t, x \in \mathbb{R}.$$

From Proposition 2.2 we obtain for $n, m \in \mathbb{N}_0$ and $\varphi, \xi \in \mathcal{N}$

$$\int_{\mathcal{N}'} H_{n, (\varphi, \varphi)_{\mathcal{H}}} (\langle \varphi, \omega \rangle) H_{m, (\xi, \xi)_{\mathcal{H}}} (\langle \xi, \omega \rangle) \, d\mu = \delta_{n,m} n! (\varphi, \xi)_{\mathcal{H}}^n = \delta_{n,m} n! \left(\varphi^{\otimes n}, \xi^{\otimes n}\right)_{\mathcal{H}^{\otimes n}}.$$

$$(4)$$

Let $n \in \mathbb{N}_0$ be fixed and $I$ a finite index set and $\varphi_i \in \mathcal{N}$, $\alpha_i \in \mathbb{C}$ for $i \in I$. Define the function in $L^2(\mu)$

$$\mathcal{N}' \ni \omega \mapsto \left\langle \sum_{i \in I} \alpha_i \varphi_i^{\otimes n}, :\omega^{\otimes n}: \right\rangle := \sum_{i \in I} \alpha_i H_{n, (\varphi_i, \varphi_i)_{\mathcal{H}}} (\langle \varphi_i, \omega \rangle) \in \mathbb{C}. \quad (5)$$

From (4) we obtain the Itô isometry between $L^2(\mu)$ and the symmetric tensor product $\mathcal{H}_{\mathbb{C}}^{\hat{\otimes} n}$

$$\left\| \left\langle \sum_{i \in I} \alpha_i \varphi_i^{\otimes n}, :\cdot^{\otimes n}: \right\rangle \right\|_{L^2(\mu)}^2 = n! \left\| \sum_{i \in I} \alpha_i \varphi_i^{\otimes n} \right\|_{\mathcal{H}_{\mathbb{C}}^{\otimes n}}^2. \quad (6)$$

Observe that we consider on the symmetric space $\mathcal{H}_{\mathbb{C}}^{\hat{\otimes} n}$ the scalar product $n!(\cdot, \cdot)_{\mathcal{H}_{\mathbb{C}}^{\otimes n}}$, where $(\cdot, \cdot)_{\mathcal{H}_{\mathbb{C}}^{\otimes n}}$ is the usual scalar product on $\mathcal{H}_{\mathbb{C}}^{\otimes n}$, see also [15, Appendix 2]. From the polarization identity we obtain that elements $\sum_{i \in I} \alpha_i \varphi_i^{\otimes n} \in \mathcal{H}_{\mathbb{C}}^{\hat{\otimes} n}$ as above form a dense subset of the complex symmetric tensor product $\mathcal{H}_{\mathbb{C}}^{\hat{\otimes} n}$. Hence, via (6) and an approximating sequence, for an element $f^{(n)} \in \mathcal{H}_{\mathbb{C}}^{\hat{\otimes} n}$ we obtain an element $F_n \in W_n(\mu)$ which we denote by $F_n = \langle f^{(n)}, :\cdot^{\otimes n}: \rangle$ satisfying

$$\|F_n\|_{L^2(\mu)}^2 = n! \left\| f^{(n)} \right\|_{\mathcal{H}_{\mathbb{C}}^{\otimes n}}^2.$$

Conversely, representing usual monomials via Hermite polynomials, we obtain that every element $F_n \in W_{(n)}(\mu)$ has a representation as $F_n = \langle f^{(n)}, :\cdot^{\otimes n}: \rangle$ where $\langle f^{(n)}, :\cdot^{\otimes n}: \rangle$ denotes again the $L^2(\mu)$-limit of elements as in (5). Let $\Gamma(\mathcal{H})$ be the *symmetric Fock space* over $\mathcal{H}$, i.e.,

$$\Gamma(\mathcal{H}) := \left\{ \mathbf{f} = (f^{(0)}, f^{(1)}, f^{(2)}, \dots) \mid f^{(n)} \in \mathcal{H}_{\mathbb{C}}^{\hat{\otimes} n} \text{ for all} \right.$$
$$\left. n \in \mathbb{N}_0, \sum_{n=0}^{\infty} n! \left\| f^{(n)} \right\|_{\mathcal{H}}^2 < \infty \right\}. \quad (7)$$

Observe that we used the abbreviation $\| f^{(n)} \|_{\mathcal{H}}$ for the norm $\| f^{(n)} \|_{\mathcal{H}_{\mathbb{C}}^{\otimes n}}$ in (7) and use henceforth similar notation for corresponding scalar products to keep the notation simple. The space $\Gamma(\mathcal{H})$ carries the scalar product

$$(\mathbf{f}, \mathbf{g})_{\Gamma(\mathcal{H})} = \sum_{n=0}^{\infty} n! (f_n, g_n)_{\mathcal{H}}, \text{ for } \mathbf{f}, \mathbf{g} \in \Gamma(\mathcal{H}).$$

The above derived decomposition of $L^2(\mu)$ is the subject of the *Wiener–Itô–Segal* theorem.

**Theorem 2.3** (Wiener–Itô–Segal isomorphism) *The mapping*

$$I : \Gamma(\mathcal{H}) \to L^2(\mu), \ \mathbf{f} \mapsto \sum_{n=0}^{\infty} \left\langle f^{(n)}, : \cdot^{\otimes n} : \right\rangle \tag{8}$$

*is a unitary isomorphism.*

Hence, each $F \in L^2(\mu)$ has a unique *chaos decomposition* $F = \sum_{n=0}^{\infty} \left\langle f^{(n)}, : \cdot^{\otimes n} : \right\rangle$ with *kernels* $f^{(n)} \in \mathcal{H}_{\mathbb{C}}^{\widehat{\otimes} n}$, $n \in \mathbb{N}_0$, and $\|F\|_{L^2(\mu)}^2 = \sum_{n=0}^{\infty} n! \left\| f^{(n)} \right\|_{\mathcal{H}}^2$. From this point we can easily define spaces of random and generalised random variables via the concepts of second quantisation, for more details see [15, Chapter 3.C]. Let $(A, D(A))$ be a closed and densely defined linear operator on $\mathcal{H}$ with $\|Af\|_{\mathcal{H}} \geq \|f\|_{\mathcal{H}}$ for all $f \in D(A)$. We define the Hilbert space $(\mathcal{G}_A, \|\cdot\|_{\mathcal{G}_A})$ as the domain of the second quantisation of $A$, i.e.,

$$\mathcal{G}_A = \left\{ F = \sum_{n=0}^{\infty} \left\langle f^{(n)}, : \cdot^{\otimes n} : \right\rangle \in L^2(\mu) \ \middle| \right.$$

$$\left. f^{(n)} \in D(A)^{\otimes n}, \sum_{n=0}^{\infty} n! \|A^{\otimes n} f^{(n)}\|_{\mathcal{H}}^2 < \infty \right\},$$

$$\|F\|_{\mathcal{G}_A}^2 := \sum_{n=0}^{\infty} n! \|A^{\otimes n} f^{(n)}\|_{\mathcal{H}}^2, \quad F \in \mathcal{G}_A.$$

The main objective of this paper is to study and characterize the space $\mathcal{G}_A$ for a special choice of $A$. To this end we first lift the rigging in (1). Therefore we need the following lemma.

**Lemma 2.4** *Let $(A_1, D(A_1))$ and $(A_2, D(A_2))$ be two closed and densely defined linear operators on $\mathcal{H}$ satisfying $\|A_i f\|_{\mathcal{H}} \geq \|f\|_{\mathcal{H}}$ for all $f \in D(A_i)$ and $i = 1, 2$. Assume that $D(A_1)$ is continuously and densely embedded into $D(A_2)$, where both spaces are equipped with $\|A_1 \cdot\|_{\mathcal{H}}$ and $\|A_2 \cdot\|_{\mathcal{H}}$, respectively. Then the space $\mathcal{G}_{A_1}$ is densely and continuously embedded into $\mathcal{G}_{A_2}$.*

**Proof** This follows as in [15, Chapter 3.B, pp. 54–55]. $\qquad\qquad\qquad\qquad\qquad \square$

From the theory of closed and symmetric bilinear forms, see e.g. [26], there exists for every $p \in \mathbb{N}$ a linear closed and densely defined linear operator $(A_p, D(A_p))$ on $\mathcal{H}$ s.t. for all $f, g \in \mathcal{H}_p$ it holds $(f, g)_p = (A_p f, A_p g)_{\mathcal{H}}$. By the previous considerations and Lemma 2.4 we can form the Hilbert spaces $(\mathcal{H}_p) := \mathcal{G}_{A_p}$ together with their dual spaces $(\mathcal{H}_{-p}) := (\mathcal{H}_p)'$, $p \in \mathbb{N}$, and obtain the chain of continuous and dense embeddings

$$(\mathcal{N}) \subseteq (\mathcal{H}_p) \subseteq (\mathcal{H}_q) \subseteq L^2(\mu) \subseteq (\mathcal{H}_{-q}) \subseteq (\mathcal{H}_{-p}) \subseteq (\mathcal{N})', \quad p \geq q,$$

where $(\mathcal{N}) = \bigcap_{p \in \mathbb{N}} (\mathcal{H}_p)$ is equipped with the projective limit topology of the spaces $\left( (\mathcal{H}_p) \right)_{p \in \mathbb{N}}$, see also [28, Sect. II.5.], and $(\mathcal{N})' = \bigcup_{p \in \mathbb{N}} (\mathcal{H}_{-p})$ is the dual space of $(\mathcal{N})$ carrying the inductive limit topology of the spaces $\left((\mathcal{H}_{-p})\right)_{p \in \mathbb{N}}$, see also [28, Sect. II.6.]. The dual pairing between elements $F \in (\mathcal{N})$ and $\Phi \in (\mathcal{N})'$ is denoted by $\langle\!\langle F, \Phi \rangle\!\rangle := \Phi(F)$.

We now specify the assumption on an operator $(K, D(K))$ for which we want to study the space $\mathcal{G}_K$ and its dual space in greater detail.

**Assumption 2.5** Let $(K, D(K))$ be a densely defined self-adjoint operator on $\mathcal{H}$ with the properties:

 (i) The spectrum $\mathrm{spec}(K)$ of $K$ satisfies $\mathrm{spec}(K) \subseteq [1, \infty)$,
 (ii) $\mathcal{N} \subseteq D(K)$ and the closure of $(K, \mathcal{N})$ equals $(K, D(K))$, i.e., $\mathcal{N}$ is a core for $(K, D(K))$,
 (iii) $K : (\mathcal{N}, d) \longrightarrow (\mathcal{N}, d)$ is continuous and bijective, where $d$ is defined in the beginning of Sect. 2.

For an operator $(K, D(K))$ satisfying Assumption 2.5 we denote by $(K^s, D(K^s))$, $s \in \mathbb{N}$, the closure of $(K^s, \mathcal{N})$ defined on $\mathcal{H}$. In particular, $(K^s, D(K^s))$ satisfies Assumption 2.5, too. Thus, by the same arguments as above, for such an operator $(K, D(K))$ we define $\mathcal{G}_{K,s} := \mathcal{G}_{K^s}$ and $\mathcal{G}_{K,-s} := \mathcal{G}'_{K,s}$. Once more, we obtain the continuous and dense embeddings

$$\mathcal{G}_K \subseteq \mathcal{G}_{K,s} \subseteq \mathcal{G}_{K,l} \subseteq L^2(\mu) \subseteq \mathcal{G}_{K,-l} \subseteq \mathcal{G}_{K,-s} \subseteq \mathcal{G}'_K, \quad s \geq l,$$

where $\mathcal{G}_K = \bigcap_{s \in \mathbb{N}} \mathcal{G}_{K,s}$ is equipped with the projective limit topology of the spaces $\left(\mathcal{G}_{K,s}\right)_{s \in \mathbb{N}}$ and $\mathcal{G}'_K = \bigcup_{s \in \mathbb{N}} \mathcal{G}_{K,-s}$ is the dual space of $\mathcal{G}_K$ carrying the inductive limit topology of the spaces $\left(\mathcal{G}_{K,-s}\right)_{s \in \mathbb{N}}$.

**Definition 2.6** For $\Phi \in (\mathcal{N})'$, its *S-transform* is defined by

$$S\Phi : \mathcal{N} \to \mathbb{C}, \varphi \mapsto \langle\!\langle :\exp\left(\langle \varphi, \cdot \rangle\right):, \Phi \rangle\!\rangle,$$

where $:\exp\left(\langle \varphi, \cdot \rangle\right): = \sum_{n=0}^{\infty} \frac{1}{n!} \left\langle \varphi^{\otimes n}, :\cdot^{\otimes n}: \right\rangle = \exp\left(\langle \varphi, \cdot \rangle - \frac{1}{2} \langle \varphi, \varphi \rangle\right) \in (\mathcal{N})$ is the Wick exponential of $\varphi \in \mathcal{N}$.

For the next theorem, we need the notion of U-functionals:

**Definition 2.7** A map $U : \mathcal{N} \to \mathbb{C}$ is called a *U-functional*, if the following two conditions are fulfilled

 (i) $U$ is *ray-entire*, i.e. for all $\varphi, \psi \in \mathcal{N}$, the function

$$\mathbb{R} \ni x \mapsto U(\varphi + x\psi)$$

extends to an entire function on $\mathbb{C}$,
 (ii) $U$ is *uniformly bounded of exponential order 2*, i.e. there exist $0 \leq A, B < \infty$ and $p \in \mathbb{N}$ s.t. for all $\varphi \in \mathcal{N}$ and $\lambda \in \mathbb{C}$ it holds

$$|U(\lambda\varphi)| \leq A \exp\left(B |\lambda|^2 \|\varphi\|_p^2\right),$$

here $U$ denotes the extension of $U$ to $\mathcal{N}_{\mathbb{C}}$ given in (i).

The following important characterisation theorem shows that there is a bijection between $(\mathcal{N})'$ and the set of U-functionals. For a proof see [18, Theorem 11].

**Theorem 2.8** *The S-transform is a bijection between* $(\mathcal{N})'$ *and the set of* $U$*-functionals.*

Our goal is to characterize the spaces $\mathcal{G}_{K,s}$, $s \in \mathbb{Z}$, in terms of $U$-functionals. To this end we first explain how the pairs of spaces $((\mathcal{N}), (\mathcal{N})')$ and $(\mathcal{G}_K, \mathcal{G}'_K)$ are related, using Lemma 2.4.

**Lemma 2.9** *Assume* $(K, D(K)$ *satisfies Assumption* 2.5*. The space* $(\mathcal{N})$ *is continuously and densely embedded into* $\mathcal{G}_K$*. Hence the following chain of continuous and dense embeddings holds true*

$$(\mathcal{N}) \subseteq \mathcal{G}_K \subseteq L^2(\mu) \subseteq \mathcal{G}'_K \subseteq (\mathcal{N})'. \tag{9}$$

*Proof* Let $s \in \mathbb{N}$ be arbitrary. Since $K^s : (\mathcal{N}, d) \longrightarrow (\mathcal{N}, d)$ is continuous there exists a $p \in \mathbb{N}$ and $C \in (0, \infty)$ s.t. $\|K^s \varphi\|_{\mathcal{H}} \leq K \|A_p \varphi\|_{\mathcal{H}}$ for all $\varphi \in \mathcal{N}$. By the closedness of $(K^s, D(K^s))$ and $(A_p, \mathcal{H}_p)$ we obtain that the norms $\|K^s \cdot\|_{\mathcal{H}}$, $\|A_p \cdot\|_{\mathcal{H}}$ are compatible on $\mathcal{N}$. Thus we obtain that $(A_1, D(A_1)) = (A_p, \mathcal{H}_p)$ and $(A_2, D(A_2)) = (K^s, D(K^s))$ satisfy the assumption of Lemma 2.4. Therefore we obtain the dense and continuous embedding of $(\mathcal{H}_p)$ into $\mathcal{G}_{K,s}$. The first embedding in (9) follows now by the definition of $(\mathcal{N})$ and $\mathcal{G}_K$. The second embedding follows by the same argument for $(K, D(K))$ and the identity operator on $\mathcal{H}$. The remaining assertions follow immediately by the previous ones. $\qquad \square$

Observe that the triple $(\mathcal{N}, \mathcal{H}, (K, D(K)))$ determines our probabilistic set up (9) completely.

**Definition 2.10** Let $m \in \mathbb{N}$ and $(\varphi_i)_{i=1}^m \subset \mathcal{N}$ be an orthonormal system in $\mathcal{H}$. We call

$$P : \mathcal{N}'_{\mathbb{C}} \to \mathcal{N}_{\mathbb{C}}, \, P\eta := \sum_{i=1}^m \langle \varphi_i, \eta \rangle \, \varphi_i$$

an *orthogonal projection* from $\mathcal{N}'_{\mathbb{C}}$ into $\mathcal{N}_{\mathbb{C}}$. We denote the set of all orthogonal projections from $\mathcal{N}'_{\mathbb{C}}$ into $\mathcal{N}_{\mathbb{C}}$ by $\mathbb{P}$.

Recall the measure $\mu_{\frac{1}{2}}$ on $\mathcal{N}'$ from above. On the complexification $\mathcal{N}'_{\mathbb{C}} = \mathcal{N}' \times \mathcal{N}'$ we define the product measure $\nu = \mu_{\frac{1}{2}} \otimes \mu_{\frac{1}{2}}$. Now we can formulate the following characterisation of the spaces $\mathcal{G}_K$ and $\mathcal{G}'_K$ which is the main result of this paper.

**Theorem 2.11** *Let* $(K, D(K))$ *satisfy Assumption* 2.5 *and* $\Phi \in (\mathcal{N})'$ *or equivalently let* $U$ *be a* $U$*-functional s.t.* $S^{-1}U = \Phi$*. Then for* $s \in \mathbb{Z}$ *the two statements*

(i) $\Phi \in \mathcal{G}_{K,s}$,
(ii) $\sup_{P \in \mathbb{P}} \int_{\mathcal{N}'_{\mathbb{C}}} |U(K^s P \eta)|^2 \, \nu(d\eta) < \infty$

*are equivalent. In particular, the following two equivalencies are true.*

(i) $\Phi \in \mathcal{G}_K \iff \forall s \in \mathbb{N} : \sup_{P \in \mathbb{P}} \int_{\mathcal{N}'_{\mathbb{C}}} |U(K^s P \eta)|^2 \, \nu(d\eta) < \infty.$

(ii) $\Phi \in \mathcal{G}'_K \Longleftrightarrow \exists s \in \mathbb{N} : \sup_{P \in \mathbb{P}} \int_{\mathcal{N}'_{\mathbb{C}}} \left| U(K^{-s} P\eta) \right|^2 \nu(d\eta) < \infty.$

*In particular, for $K = Id$ and $\Phi \in (\mathcal{N})'$ we obtain the following equivalence*

$$\Phi \in L^2(\mu) \Longleftrightarrow \sup_{P \in \mathbb{P}} \int_{\mathcal{N}'_{\mathbb{C}}} |U(P\eta)|^2 \nu(d\eta) < \infty.$$

Before proceeding we present some classical examples of the functional analytic framework $(\mathcal{H}, \mathcal{N})$ as well as interesting choices for $(K, D(K))$. In particular, in Example 2.12(ii) we show how to construct the necessary nuclear rigging (1) for typical examples arising in the contexts of stochastic partial differential equations, see Sect. 5.2.

**Example 2.12** (i) Let $d_1, d_2 \in \mathbb{N}$. The real Hilbert space $\mathcal{H} := L^2(\mathbb{R}^{d_1}; \mathbb{R}^{d_2})$ and the nuclear space $\mathcal{N} = S(\mathbb{R}^{d_1}; \mathbb{R}^{d_2})$ of square integrable functions and Schwartz test functions mapping from $\mathbb{R}^{d_1}$ to $\mathbb{R}^{d_2}$ respectively. Different examples for a family of seminorms $(\|\cdot\|_p)_{p \in \mathbb{N}_0}$ defined on $S(\mathbb{R}^{d_1}; \mathbb{R}^{d_2})$ satisfying the assumptions above can be found in [26, Appendix to V.3]. For $d_1 = d_2 = 1$ this setting is called the standard White noise setting. Observe that the process $\left( \left\langle \left( \otimes_{i=1}^{d_1} 1_{[0,t_i]} \right)^{\times d_2}, \cdot \right\rangle \right)_{t_1, \ldots, t_{d_1} \geq 0} \subseteq L^2(\mathcal{N}', \mu)$ is a modification of a $(d_1, d_2)$-Brownian sheet.

(ii) Let $d_1, d_2 \in \mathbb{N}$. More generally as in (i) assume that $\sigma$ is a tempered measure on $\mathbb{R}^{d_1}$, i.e., $\int_{\mathbb{R}^{d_1}} \frac{1}{1+|\xi|^{2m}} \sigma(d\xi) < \infty$ for some $m \in \mathbb{N}$, which has full topological support on $\mathbb{R}^d$. Define $\mathcal{H}$ to be the completion of $\mathcal{N} := S(\mathbb{R}^{d_1}; \mathbb{R}^{d_2})$ w.r.t. the norm $\|\varphi\|_{\mathcal{H}} := \left( \int_{\mathbb{R}^{d_1}} |\mathcal{F}\varphi|^2 \, d\sigma \right)^{\frac{1}{2}}$, where $\mathcal{F}$ denotes the component-wise Fourier transform and $|\cdot|$ the euclidean norm on $\mathbb{R}^{d_2}$. Let $(\|\cdot\|_p)_{p \in \mathbb{N}}$ denote an increasing family of consistent seminorms on $S(\mathbb{R}^{d_1}; \mathbb{R}^{d_2})$ inducing the Schwartz space topology. Since $\sigma$ is tempered there exists a $q \in \mathbb{N}$ s.t. $\max\{\|\varphi\|_{L^1(\mathbb{R}^{d_1}, \mathbb{R}^{d_2})}, \|\varphi\|_{\mathcal{H}}\} \leq \|\varphi\|_q$ for all $\varphi \in S(\mathbb{R}^{d_1}, \mathbb{R}^{d_2})$. From the full support of $\sigma$ and the continuity of $\mathcal{F} : L^1(\mathbb{R}^{d_1}, \mathbb{R}^{d_2}) \longrightarrow L^\infty(\mathbb{R}^{d_1}, \mathbb{R}^{d_2})$ one concludes the consistency of the norms $\|\cdot\|_q$ and $\|\cdot\|_{\mathcal{H}}$ on $S(\mathbb{R}^{d_1}; \mathbb{R}^{d_2})$. Thus in the sense of Gelfand triples we obtain the dense embeddings

$$S(\mathbb{R}^{d_1}; \mathbb{R}^{d_2}) \subseteq \mathcal{H}_q \subseteq \mathcal{H} \subseteq \mathcal{H}_{-q} \subseteq S'(\mathbb{R}^{d_1}; \mathbb{R}^{d_2}).$$

Consequently we obtain the Gaussian measure $\mu$ on $S'(\mathbb{R}^{d_1}; \mathbb{R}^{d_2})$ with covariance $(\cdot, \cdot)_{\mathcal{H}}$, i.e., $\mu$ satisfies

$$\int_{S'(\mathbb{R}^{d+1})} \exp(i \langle \varphi, \omega \rangle) \mu(d\omega) = \exp(-\frac{1}{2}(\varphi, \varphi)_{\mathcal{H}}), \quad \varphi \in S(\mathbb{R}^{d_1}; \mathbb{R}^{d_2}).$$

In particular, the elements from the first order chaos $\langle h, \cdot \rangle \in W_{(1)} \subseteq L^2(S'(\mathbb{R}^{d_1}; \mathbb{R}^{d_2}), \mu)$, $h \in \mathcal{H}$, define a Gaussian process index by $\mathcal{H}$. In this way a huge variety of Gaussian processes, such as fractional Brownian motion, can be constructed.

(iii) An important choice for the operator $K$ is given by a multiple $\lambda > 1$ of the identity operator $Id$ on $\mathcal{H}$, i.e., $K = \lambda Id$. The space $\mathcal{G}_{\lambda Id}$ was systematically introduced in the White Noise setting in [25]. An important feature of $\mathcal{G}_{\lambda Id}$ is that this space is densely and continuously embedded into the Meyer-Watanabe space $\mathbb{D}$, see the last mentioned reference. Thus, elements from $\mathcal{G}_{\lambda Id}$ are infinitely often Malliavin differentiable and the Malliavin derivatives of arbitrary order are contained in $L^p(\mathcal{N}', \mu)$ for every $p \in [1, \infty)$.

## 3 Generalised chaos decomposition and Gaussian analysis on complex spaces

In this section we state some additional aspects of Gaussian Analysis. For further reading, see e.g. [2,15,19,22].

### 3.1 Generalised chaos decomposition

Next we generalise the chaos decomposition (8) of elements from $L^2(\mu)$ to elements from the dual spaces $\mathcal{G}_{K,-s}$, $s \in \mathbb{N}$. Let $s \in \mathbb{N}$ and recall that $\mathcal{G}_{K,s}$ is isometrically isomorphic to $\Gamma(D(K^s))$. Hence, the dual space $\mathcal{G}'_{K,s} = \mathcal{G}_{K,-s}$ is isometrically isomorphic to $(\Gamma(D(K^s)))' \cong \Gamma(D(K^s)')$, where $\cong$ denotes an isometric isomorphism and $D(K^s)'$ is the dual space of $(D(K^s), (K^s\cdot, K^s\cdot)_{\mathcal{H}})$. Observe that $D(K^s)'$ is isometrically isomorphic to the completion of $\mathcal{H}$ w.r.t. the inner product $(K^{-s}\cdot, K^{-s}\cdot)_{\mathcal{H}}$. An element $\Phi \in \mathcal{G}_{K,-s}$ which is in correspondence with $(\Phi^{(n)})_{n \in \mathbb{N}_0} \in \Gamma(D(K^s)')$ we also denote by

$$\Phi = \sum_{n=0}^{\infty} \left\langle \Phi^{(n)}, :\cdot^{\otimes n}: \right\rangle. \tag{10}$$

The correspondence in (10) is called the *generalised chaos decomposition* of $\Phi$. The dual pairing between $\Phi$ and an element $\psi \in \mathcal{G}_{K,s}$ with chaos decomposition $\psi = \sum_{n=0}^{\infty} \left\langle \psi^{(n)}, :\cdot^{\otimes n}: \right\rangle$ is given by

$$\langle\!\langle \psi, \Phi \rangle\!\rangle = \sum_{n=0}^{\infty} n! \left\langle \psi^{(n)}, \Phi^{(n)} \right\rangle, \tag{11}$$

where the dual pairing $\left\langle \psi^{(n)}, \Phi^{(n)} \right\rangle$ on the right-hand side of (11) is the one between the Hilbert space $\left( D((K^s)^{\otimes n}), \|(K^s)^{\otimes n}\cdot\|_{\mathcal{H}} \right)$ and its dual space $D\left( (K^s)^{\otimes n} \right)'$, for $n \in \mathbb{N}_0$.

From Lemma 2.9 we obtain that an element $\Phi \in \mathcal{G}'_K$ has a well-defined $S$-transform $S\Phi : \mathcal{N} \longrightarrow \mathbb{C}$. If $\Phi \in \mathcal{G}'_K$ has the generalised chaos decomposition $\Phi = \sum_{n=0}^{\infty} \left\langle \Phi^{(n)}, :\cdot^{\otimes n}: \right\rangle$ then the $S$-transform is given by

$$S\Phi(\varphi) = \sum_{n=0}^{\infty} \left\langle \varphi^{\otimes n}, \Phi^{(n)} \right\rangle, \quad \varphi \in \mathcal{N}_{\mathbb{C}}. \tag{12}$$

### 3.2 Gaussian analysis on complex spaces

In this part we briefly present the analogue of the orthogonal decomposition of $L^2(\mu)$ for a closed subspace $E^2(\nu)$ of $L^2(\mathcal{N}'_{\mathbb{C}}, \nu)$. The major difference between the space $L^2(\mu)$ and $E^2(\nu)$ is that in the latter case there is no need for using Hermite polynomials, see Proposition 3.2. The underlying reason is that the monomials of different order automatically form an orthogonal system in $L^2(\mathbb{C}, e^{-|z|^2_{euc}} dz)$. The proofs of the next two propositions are elementary and therefore we skip them.

**Proposition 3.1** *Let $\varphi_1, ..., \varphi_n \in \mathcal{N}$, $n \in \mathbb{N}$. The image measure of $\nu$ under the map*

$$T_{\varphi_1, ..., \varphi_n} : \mathcal{N}'_{\mathbb{C}} \longrightarrow \mathbb{C}^n, \eta \mapsto (\langle \varphi_i, \eta \rangle)_{i=1, ..., n}$$

*is absolutely continuous w.r.t. the Lebesgue measure $dz$ on $\mathbb{C}^n$ and has the Radon-Nikodym derivative*

$$\frac{d\nu \circ T_{\varphi_1, ..., \varphi_n}^{-1}}{dz}(z) = \frac{1}{\pi^n} e^{-\bar{z}^T C z}, \quad z \in \mathbb{C}^n$$

*where $C = \left( (\varphi_i, \varphi_j)_{\mathcal{H}} \right)_{1 \le i, j \le n} \in \mathbb{R}^{n \times n}$.*

The space of polynomials $\mathcal{P}(\mathcal{N}'_{\mathbb{C}})$ on $\mathcal{N}'_{\mathbb{C}}$ is given by collection of all functions $G : \mathcal{N}'_{\mathbb{C}} \to \mathbb{C}$ which are given as $G(\eta) = p(\langle \varphi_1, \eta \rangle, ..., \langle \varphi_k, \eta \rangle)$, where $p$ is a complex polynomial in $k \in \mathbb{N}$ variables and $\varphi_i \in \mathcal{N}$, for $i = 1, .., k$.

**Proposition 3.2** *Let $m, n \in \mathbb{N}$, $\varphi, \psi \in \mathcal{N}$. Then it holds*

$$\left( \langle \varphi, \cdot \rangle^n, \langle \psi, \cdot \rangle^m \right)_{L^2(\nu)} = \delta_{m,n} \cdot n! \cdot (\varphi^{\otimes n}, \psi^{\otimes n})_{\mathcal{H}}. \tag{13}$$

*In particular, $\mathcal{P}(\mathcal{N}'_{\mathbb{C}}) \subseteq L^2(\nu)$.*

Similar as in the derivation of Theorem 2.3, for $f^{(n)} \in \mathcal{H}_{\mathbb{C}}^{\widehat{\otimes} n}$ we can define an element in $L^2(\nu)$ denoted by $\left\langle f^{(n)}, \cdot^{\otimes n} \right\rangle$ which is given as the $L^2(\nu)$-limit of polynomials, i.e.,

$$\left\langle f^{(n)}, \cdot^{\otimes n} \right\rangle := \lim_{m \to \infty} \sum_{k=1}^{l_m} \alpha_{k,m} \left\langle \varphi_{k,m}, \cdot \right\rangle^n \in L^2(\nu),$$

where $l_m \in \mathbb{N}$, $\alpha_{k,m} \in \mathbb{C}$, $\varphi_{k,m} \in \mathcal{N}$ for all $k = 1, ..., l_m$, $m \in \mathbb{N}$, and it holds

$$f^{(n)} = \lim_{m \to \infty} \sum_{k=1}^{l_m} \alpha_{k,m} \varphi_{k,m}^{\otimes n} \in \mathcal{H}_{\mathbb{C}}^{\widehat{\otimes} n}.$$

In particular, the orthogonality relation (13) stays valid in the limit case, i.e., for $f^{(n)}, g^{(n)} \in \mathcal{H}_{\mathbb{C}}^{\widehat{\otimes} n}$ it holds

$$\left( \left\langle f^{(n)}, \cdot^{\otimes n} \right\rangle, \left\langle g^{(n)}, \cdot^{\otimes n} \right\rangle \right)_{L^2(\nu)} = \delta_{m,n} \cdot n! \cdot (f^{(n)}, g^{(n)})_{\mathcal{H}}. \tag{14}$$

In contrast to the real case, the polynomials $\mathcal{P}(\mathcal{N}'_{\mathbb{C}})$ are not dense in $L^2(\nu)$. Their closure is the so called Bargmann-Segal space $E^2(\nu)$, see also [12], which is given by

$$E^2(\nu) := \overline{\mathcal{P}(\mathcal{N}'_{\mathbb{C}})}^{L^2(\nu)} = \Big\{ \sum_{n=0}^{\infty} \left\langle g^{(n)}, \cdot^{\otimes n} \right\rangle \mid g^{(n)} \in \mathcal{H}_{\mathbb{C}}^{\widehat{\otimes} n}, \ \sum_{n=0}^{\infty} n! \left\| g^{(n)} \right\|_{\mathcal{H}}^2 < \infty \Big\}.$$

## 4 Proof of Theorem 2.11

This section is devoted to the proof of Theorem 2.11, which is our main result.

Recall the chain of continuous embeddings from (1). This chain lifts to the $n$-fold symmetric complexified tensor powers, see e.g. [15, Chapter 3.B], i.e, we obtain continuous embeddings

$$\mathcal{N}_{\mathbb{C}}^{\widehat{\otimes} n} \subseteq \mathcal{H}_{p,\mathbb{C}}^{\widehat{\otimes} n} \subseteq \mathcal{H}_{q,\mathbb{C}}^{\widehat{\otimes} n} \subseteq \mathcal{H}_{\mathbb{C}}^{\widehat{\otimes} n} \subseteq \mathcal{H}_{-q,\mathbb{C}}^{\widehat{\otimes} n} \subseteq \mathcal{H}_{-p,\mathbb{C}}^{\widehat{\otimes} n} \subseteq \mathcal{N}_{\mathbb{C}}'^{\widehat{\otimes} n}, \quad p \geq q,$$

where $\mathcal{N}_{\mathbb{C}}^{\widehat{\otimes} n} := \bigcap_{p \in \mathbb{N}} \mathcal{H}_{p,\mathbb{C}}^{\widehat{\otimes} n}$ is equipped with the projective limit topology of the Hilbert spaces $\mathcal{H}_{p,\mathbb{C}}^{\widehat{\otimes} n}$, $p \in \mathbb{N}$ and $\mathcal{N}_{\mathbb{C}}'^{\widehat{\otimes} n}$ is the dual space of $\mathcal{N}_{\mathbb{C}}^{\widehat{\otimes} n}$ which satisfies $\mathcal{N}_{\mathbb{C}}'^{\widehat{\otimes} n} = \bigcup_{p \in \mathbb{N}} \mathcal{H}_{-p,\mathbb{C}}^{\widehat{\otimes} n}$ and carries the inductive limit topology of the spaces $\mathcal{H}_{-p,\mathbb{C}}^{\widehat{\otimes} n}$, $p \in \mathbb{N}$.

The operator $K : \mathcal{N} \longrightarrow \mathcal{N}$ was assumed to be bijective and continuous, hence by the inverse mapping theorem $K^s$, $s \in \mathbb{Z}$, is also continuous, see [27, Corollary I.2.12(b)]. By the same procedure which leads to tensor powers of operators between Hilbert spaces, we can define $(K^s)^{\otimes n}$ for $s \in \mathbb{Z}$ and $n \in \mathbb{N}$ as a well-defined, linear and continuous operator on $\mathcal{N}_{\mathbb{C}}^{\widehat{\otimes} n}$. Observe that $(K^s)^{\otimes n}$ is bijective from $\mathcal{N}_{\mathbb{C}}^{\widehat{\otimes} n}$ into itself and the tensor powers $\left( (K^s)^{\otimes n}, D((K^s)^{\otimes n}) \right)$ are self-adjoint on $\mathcal{H}_{\mathbb{C}}^{\widehat{\otimes} n}$, where $D((K^s)^{\otimes n}) = \mathcal{H}_{\mathbb{C}}^{\widehat{\otimes} n}$ if $s \leq 0$, for all $s \in \mathbb{Z}$. Hence, for all $s \in \mathbb{Z}$ and $n \in \mathbb{N}$ we can define an extension of $(K^s)^{\otimes n}$ to $\mathcal{N}_{\mathbb{C}}'^{\widehat{\otimes} n}$ in the following way:

$$(K^s)^{\otimes n} : \mathcal{N}_{\mathbb{C}}'^{\widehat{\otimes} n} \longrightarrow \mathcal{N}_{\mathbb{C}}'^{\widehat{\otimes} n}, \Phi \mapsto (K^s)^{\otimes n} \Phi := \Phi \circ (K^s)^{\otimes n}. \tag{15}$$

For the next proposition recall that every element $\Phi \in (\mathcal{N})'$ has a generalised chaos decomposition $\Phi = \sum_{n=0}^{\infty} \langle \Phi^{(n)}, : \cdot^{\otimes n} : \rangle$ where for some $p \in \mathbb{N}$ it holds $\Phi^{(n)} \in \mathcal{H}_{-p,\mathbb{C}}^{\widehat{\otimes} n} \subseteq \mathcal{N}_{\mathbb{C}}'^{\widehat{\otimes} n}$ for all $n \in \mathbb{N}_0$.

**Proposition 4.1** *Let $s \in \mathbb{Z}$. Then it holds*

$$\mathcal{G}_{K,s} = \left\{ \Phi = \sum_{n=0}^{\infty} \left\langle \Phi^{(n)}, : \cdot^{\otimes n} : \right\rangle \in (\mathcal{N})' \; \middle| \right.$$

$$\left. (K^s)^{\otimes n} \Phi^{(n)} \in \mathcal{H}_{\mathbb{C}}^{\widehat{\otimes} n}, \sum_{n=0}^{\infty} n! \left\| (K^s)^{\otimes n} \Phi^{(n)} \right\|_{\mathcal{H}}^2 < \infty \right\}, \tag{16}$$

*where $(K^s)^{\otimes n} \Phi^{(n)}$ in* (16) *is defined via* (15).

***Proof*** Denote the set on the right-hand side of (16) by $\mathcal{A}_s$. We split the proof into two parts. First let $s$ be non-negative. In this case the inclusion $\mathcal{G}_{K,s} \subseteq \mathcal{A}_s$ follows immediately by the definition of $\mathcal{G}_{K,s}$. Now let $\Phi \in \mathcal{A}_s$, i.e., $\Phi = \sum_{n=0}^{\infty} \langle \Phi^{(n)}, : \cdot^{\otimes n} : \rangle \in (\mathcal{N})'$ s.t. $(K^s)^{\otimes n} \Phi^{(n)} \in \mathcal{H}_{\mathbb{C}}^{\widehat{\otimes} n}$ and $\sum_{n=0}^{\infty} n! \left\| (K^s)^{\otimes n} \Phi^{(n)} \right\|_{\mathcal{H}}^2 < \infty$. To prove that $\Phi \in \mathcal{G}_{K,s}$ it suffices to show that $\Phi^{(n)} \in D(K^s)^{\widehat{\otimes} n}$ for all $n \in \mathbb{N}$.

By assumption, for all $n \in \mathbb{N}$ there exists a $\psi^{(n)} \in \mathcal{H}_{\mathbb{C}}^{\widehat{\otimes} n}$ s.t.

$$\left( \varphi^{(n)}, \psi^{(n)} \right)_{\mathcal{H}} = \left\langle \varphi^{(n)}, (K^s)^{\otimes n} \Phi^{(n)} \right\rangle = \left\langle (K^s)^{\otimes n} \varphi^{(n)}, \Phi^{(n)} \right\rangle, \quad \forall \varphi^{(n)} \in \mathcal{N}_{\mathbb{C}}^{\widehat{\otimes} n}.$$

Using [26, Theorem VIII.33] we obtain that $(K^s)^{\otimes n} : D((K^s)^{\otimes n}) \longrightarrow \mathcal{H}_{\mathbb{C}}^{\widehat{\otimes} n}$ is bijective and self-adjoint. Hence, we can find a $\tilde{\psi}^{(n)} \in D((K^s)^{\otimes n})$ s.t. $(K^s)^{\otimes n} \tilde{\psi}^{(n)} = \psi^{(n)}$. From the self-adjointness of $(K^s)^{\otimes n}$ we can conclude $\Phi^{(n)} = \overline{\tilde{\psi}^{(n)}}$, where $\bar{\cdot}$ is the natural complex conjugation on the complexified vector space $\mathcal{H}_{\mathbb{C}}^{\widehat{\otimes} n}$, which finishes the proof for non-negative $s$.

For the second part we replace $s$ by $-s$, $s \in \mathbb{N}$. Recall that $\mathcal{G}_{K,-s} \cong \Gamma(D(K^s)')$. Denote by $\overline{\mathcal{H}_{\mathbb{C}}^{\widehat{\otimes} n}}^{\|(K^{-s})^{\otimes n} \cdot \|}$ the abstract completion of $\mathcal{H}_{\mathbb{C}}^{\widehat{\otimes} n}$ w.r.t. $\|(K^{-s})^{\otimes n} \cdot \|$. One easily checks via the Riesz isomorphism that

$$\overline{\mathcal{H}_{\mathbb{C}}^{\widehat{\otimes} n}}^{\|(K^{-s})^{\otimes n} \cdot \|} \ni (\Phi_k^{(n)})_{k \in \mathbb{N}} \mapsto \lim_{k \to \infty} \left( (K^{-s})^{\otimes n} \cdot, (K^{-s})^{\otimes n} \Phi_k^{(n)} \right)_{\mathcal{H}} \in \left( D((K^s)^{\otimes n}) \right)' \tag{17}$$

is an isometric complex conjugate linear isomorphism. Hence, the inclusion $\mathcal{G}_{K,-s} \subseteq \mathcal{A}_{-s}$ follows. Now let $\Phi \in \mathcal{A}_{-s}$ with generalised chaos decomposition $\Phi = \sum_{n=0}^{\infty} \langle \Phi^{(n)}, : \cdot^{\otimes n} : \rangle$. It suffices to show $\Phi^{(n)} \in \left( D((K^s)^{\otimes n}) \right)'$ for all $n \in \mathbb{N}$. By assumption, for every $n \in \mathbb{N}$ there exists a $\psi^{(n)} \in \mathcal{H}_{\mathbb{C}}^{\widehat{\otimes} n}$ s.t.

$$\left\langle \varphi^{(n)}, \Phi^{(n)} \right\rangle = \left\langle (K^{-s})^{\otimes n} (K^s)^{\otimes n} \varphi^{(n)}, \Phi^{(n)} \right\rangle = \left( (K^s)^{\otimes n} \varphi^{(n)}, \psi^{(n)} \right)_{\mathcal{H}}, \quad \forall \varphi^{(n)} \in \mathcal{N}_{\mathbb{C}}^{\widehat{\otimes} n}.$$

Since $D((K^s)^{\otimes n}) = (K^{-s})^{\otimes n} \mathcal{H}_{\mathbb{C}}^{\widehat{\otimes}n}$ is dense in $\mathcal{H}_{\mathbb{C}}^{\widehat{\otimes}n}$ there exists a sequence $(\chi_k^{(n)})_{k \in \mathbb{N}}$ in $\mathcal{H}_{\mathbb{C}}^{\widehat{\otimes}n}$ s.t. $(K^{-s})^{\otimes n} \chi_k^{(n)} \longrightarrow \psi^{(n)}$ as $k \to \infty$ in $\mathcal{H}_{\mathbb{C}}^{\widehat{\otimes}n}$ for all $n \in \mathbb{N}$. Hence, by (17) we obtain $\Phi^{(n)} \in \left(D((K^s)^{\otimes n})\right)'$ which finishes the proof. $\qquad\square$

Now let $n \in \mathbb{N}$ and $P \in \mathbb{P}$ be given by $P = \sum_{j=1}^m \langle \varphi_j, \cdot \rangle \varphi_j$, where $(\varphi_j)_{j=1}^m \subset \mathcal{N}$ is an orthonormal system in $\mathcal{H}$. We consider $P$ as an orthogonal projection on $\mathcal{H}_{\mathbb{C}}$ onto the closed subspace $\mathrm{span}_{\mathbb{C}}\{\varphi_j, \, j = 1, ..., m\}$. Observe that the $n$th tensor power $P^{\otimes n}$ of $P$ defines a orthogonal projection onto the closed subspace $\mathrm{span}_{\mathbb{C}} \left\{ \hat{\otimes}_{i=1}^n \varphi_{j_i} \mid j_i \in \{1, ..., m\} \text{ for } i = 1, ..., n \right\}$ of $\mathcal{H}_{\mathbb{C}}^{\widehat{\otimes}n}$, where

$$\hat{\otimes}_{i=1}^n \varphi_{j_i} = \frac{1}{n!} \sum_{\sigma \in \mathbb{S}_n} \otimes_{i=1}^n \varphi_{j_{\sigma(i)}}$$

and $\mathbb{S}_n$ denotes the set of all permutations of $n$ elements in the following way. We extend $P^{\otimes n}$ to a linear operator on $\mathcal{N}_{\mathbb{C}}'^{\widehat{\otimes}n}$ via

$$P^{\otimes n} : \mathcal{N}_{\mathbb{C}}'^{\widehat{\otimes}n} \longrightarrow \mathcal{N}_{\mathbb{C}}^{\widehat{\otimes}n}, \Phi \mapsto \sum_{\substack{\alpha \in \{0,1,...n\}^m \\ \sum_{i=1}^m \alpha_i = n}} \left\langle \widehat{\otimes_{i=1}^m \varphi_i^{\otimes \alpha_i}}, \Phi \right\rangle \widehat{\otimes_{i=1}^m \varphi_i^{\otimes \alpha_i}}.$$

Observe that for $n \in \mathbb{N}$, $\Phi^{(n)}, \Psi^{(n)} \in \mathcal{N}_{\mathbb{C}}'^{\widehat{\otimes}n}$ and $P \in \mathbb{P}$ given as above it holds $\langle P^{\otimes n} \Phi^{(n)}, \Psi^{(n)} \rangle = \langle P^{\otimes n} \Psi^{(n)}, \Phi^{(n)} \rangle$.

**Lemma 4.2** *Let* $\Phi^{(n)} \in \mathcal{N}_{\mathbb{C}}'^{\widehat{\otimes}n}$, $n \in \mathbb{N}_0$, *fulfil* $\sup_{P \in \mathbb{P}} \sum_{n=0}^\infty n! \| P^{\otimes n} \Phi^{(n)} \|_{\mathcal{H}}^2 < \infty$. *Then it holds*

$$\left( \Phi^{(n)} \right)_{n \in \mathbb{N}_0} \in \Gamma(\mathcal{H}) \text{ and } \sup_{P \in \mathbb{P}} \sum_{n=0}^\infty n! \| P^{\otimes n} \Phi^{(n)} \|_{\mathcal{H}}^2 = \sum_{n=0}^\infty n! \| \Phi^{(n)} \|_{\mathcal{H}}^2.$$

**Proof** Since $\mathcal{N}$ is separable we can choose a dense set $\{\tilde{e}_k\}_{k \in \mathbb{N}}$ of $\mathcal{N}$. Applying the Gram-Schmidt procedure to $\{\tilde{e}_k\}_{k \in \mathbb{N}}$ we obtain a orthonormal basis $(e_k)_{k \in \mathbb{N}}$ of $\mathcal{H}$ s.t. $e_k \in \mathcal{N}$ for all $k \in \mathbb{N}$.

Define for $l \in \mathbb{N}$ the projection $P_l := \sum_{k=1}^l \langle e_k, \cdot \rangle e_k \in \mathbb{P}$. By assumption we know that the sequence $\left( (P_l^{\otimes n} \Phi^{(n)})_{n \in \mathbb{N}} \right)_{l \in \mathbb{N}} \subseteq \Gamma(\mathcal{H})$ is bounded. Therefore we can find a weakly convergent subsequence $\left( \left( P_{l_m}^{\otimes n} \Phi^{(n)} \right)_{n \in \mathbb{N}} \right)_{m \in \mathbb{N}}$ with weak limit $\left( g^{(n)} \right)_{n \in \mathbb{N}_0} \in \Gamma(\mathcal{H})$. In particular $P_{l_m}^{\otimes n} \Phi^{(n)}$ converges weakly to $g^{(n)} \in \mathcal{H}_{\mathbb{C}}^{\widehat{\otimes}n}$ as $m \to \infty$ for all $n \in \mathbb{N}_0$, i.e.,

$$\left\langle \varphi^{(n)}, P_{l_m}^{\otimes n} \Phi^{(n)} \right\rangle = \left( \varphi^{(n)}, \overline{P_{l_m}^{\otimes n} \Phi^{(n)}} \right)_{\mathcal{H}} \xrightarrow{m \to \infty} \left( \varphi^{(n)}, \overline{g^{(n)}} \right)_{\mathcal{H}} \quad \text{for all } \varphi \in \mathcal{H}_{\mathbb{C}}^{\widehat{\otimes}n}.$$

It is clear that $\Phi^{(n)}$ and $g^{(n)}$ coincide as distributions on the set
$\left\{ \widehat{\otimes_{i=1}^{\infty} e_i^{\otimes \alpha_i}} \mid \alpha \in \mathbb{N}^{\mathbb{N}}, \sum_{i=1}^{\infty} \alpha_i = n \right\}$ which is total in $\mathcal{N}_{\mathbb{C}}^{\widehat{\otimes} n}$ by the choice of $(e_k)_{k \in \mathbb{N}}$.
Thus $\Phi^{(n)} = g^{(n)} \in \mathcal{H}_{\mathbb{C}}^{\widehat{\otimes} n}$.

The last part of statement follows by the weak lower semicontinuity of the norm and the fact that for $P \in \mathbb{P}$ the restriction of $P^{\otimes n}$ to $\mathcal{H}_{\mathbb{C}}^{\widehat{\otimes} n}$ is an orthogonal projection.
□

Now we are ready to prove the main result.

**Proof of Theorem 2.11** Recall that $\mathcal{G}_K = \bigcap_{s \in \mathbb{N}} \mathcal{G}_{K,s}$ and $\mathcal{G}'_K = \bigcup_{s \in \mathbb{N}} \mathcal{G}_{K,-s}$. Hence, it suffices to show for $s \in \mathbb{Z}$ that it holds

$$\Phi \in \mathcal{G}_{K,s} \iff \sup_{P \in \mathbb{P}} \int_{\mathcal{N}'_{\mathbb{C}}} \left| S\Phi(K^s P \eta) \right|^2 \nu(d\eta) < \infty. \tag{18}$$

We make some observations which rely on (12), (14) and Lemma 4.2. Now let $s \in \mathbb{Z}$ and $\Phi = \sum_{n=0}^{\infty} \left\langle \Phi^{(n)}, :\cdot^{\otimes n}: \right\rangle \in \mathcal{G}_{K,s}$. Then it holds

$$\begin{aligned}
\|\Phi\|_{K,s}^2 &= \sum_{n \in \mathbb{N}} n! \|(K^s)^{\otimes n} \Phi^{(n)}\|_{\mathcal{H}}^2 \\
&= \sup_{P \in \mathbb{P}} \sum_{n \in \mathbb{N}} n! \|P^{\otimes n} (K^s)^{\otimes n} \Phi^{(n)}\|_{\mathcal{H}}^2 \\
&= \sup_{P \in \mathbb{P}} \int_{\mathcal{N}'_{\mathbb{C}}} \left| \sum_{n \in \mathbb{N}} \left\langle P^{\otimes n} (K^s)^{\otimes n} \Phi^{(n)}, \eta^{\otimes n} \right\rangle \right|^2 \nu(d\eta) \\
&= \sup_{P \in \mathbb{P}} \int_{\mathcal{N}'_{\mathbb{C}}} \left| \sum_{n \in \mathbb{N}} \left\langle (K^s P \eta)^{\otimes n}, \Phi^{(n)} \right\rangle \right|^2 \nu(d\eta) = \sup_{P \in \mathbb{P}} \int_{\mathcal{N}'_{\mathbb{C}}} \left| S\Phi(K^s P \eta) \right|^2 \nu(d\eta).
\end{aligned}$$

Now let $\Phi = \sum_{n=0}^{\infty} \left\langle \Phi^{(n)}, :\cdot^{\otimes n}: \right\rangle \in (\mathcal{N})'$.

The same calculations yield

$$\begin{aligned}
\sup_{P \in \mathbb{P}} \int_{\mathcal{N}'_{\mathbb{C}}} \left| S\Phi(K^s P \eta) \right|^2 \nu(d\eta) &= \sup_{P \in \mathbb{P}} \int_{\mathcal{N}'_{\mathbb{C}}} \left| \sum_{n \in \mathbb{N}} \left\langle (K^s P \eta)^{\otimes n}, \Phi^{(n)} \right\rangle \right|^2 \nu(d\eta) \\
&= \sup_{P \in \mathbb{P}} \int_{\mathcal{N}'_{\mathbb{C}}} \left| \sum_{n \in \mathbb{N}} \left\langle P^{\otimes n} (K^s)^{\otimes n} \Phi^{(n)}, \eta^{\otimes n} \right\rangle \right|^2 \nu(d\eta) \\
&= \sup_{P \in \mathbb{P}} \sum_{n \in \mathbb{N}} n! \|P^{\otimes n} (K^s)^{\otimes n} \Phi^{(n)}\|_{\mathcal{H}}^2,
\end{aligned}$$

in the second line we used the definition of $(K^s)^{\otimes n}$ given in (15). Hence, if the left-hand side is finite we obtain by Lemma 4.2 that $\left( (K^s)^{\otimes n} \Phi^{(n)} \right)_{n \in \mathbb{N}} \in \Gamma(\mathcal{H})$ which implies $\Phi \in \mathcal{G}_{K,s}$ by Proposition 4.1.
□

From the proof of Theorem 2.11 we seek the following corollary.

**Corollary 4.3** *Assume the assumptions of Theorem 2.11 are satisfied and $\Phi \in (\mathcal{N})'$. Let $(e_k)_{k \in \mathbb{N}}$ be an orthonormal basis of $\mathcal{H}$ contained in $\mathcal{N}$ which is chosen as in the proof of Lemma 4.2. Further for $l \in \mathbb{N}$ let $P_l := \sum_{k=1}^{l} \langle e_k, \cdot \rangle e_k \in \mathbb{P}$. The real sequence $\int_{\mathcal{N}_\mathbb{C}'} |S\Phi(K^s P_l \eta)|^2 \, \nu(d\eta), l \in \mathbb{N}$, is increasing in $l$. Thus, the limit as $l$ tends to infinity exists in $[0, \infty]$. Further, for $s \in \mathbb{Z}$ the statement*

$$\lim_{l \to \infty} \int_{\mathcal{N}_\mathbb{C}'} \left| S\Phi(K^s P_l \eta) \right|^2 \nu(d\eta) < \infty, \tag{19}$$

*is equivalent to the statements in (18). If (19) is satisfied we also obtain*

$$\|\Phi\|_{K,s}^2 = \lim_{l \to \infty} \int_{\mathcal{N}_\mathbb{C}'} \left| S\Phi(K^s P_l \eta) \right|^2 \nu(d\eta) < \infty.$$

**Remark 4.4** Via the spectral theorem for self-adjoint operators one could also introduce the space $\mathcal{G}_{K,s}$ for $s \in \mathbb{R}$. In the proof of Theorem 2.11 we used that the operator $K^s$, $s \in \mathbb{Z}$, maps $\mathcal{N}$ continuously into itself. If for $s \in \mathbb{R}$ the operator $K^s$ maps $\mathcal{N}$ continuously into itself the exact same proof as above also leads to the corresponding statement in Theorem 2.11 for $s \in \mathbb{R}$.

# 5 Applications

## 5.1 An equation from turbulent transport and its regularity

In this subsection we shall apply the derived characterisation to an equation from turbulent transport. Here, we specify our general setting to be the white noise setting, i.e., $\mathcal{H} = L^2(\mathbb{R})$ and $\mathcal{N} = S(\mathbb{R})$, see Example 2.12(i). Hence, we write $S'(\mathbb{R})$ for $\mathcal{N}'$, the measure $\mu$ is called *white noise measure*. We consider the case $K = \sqrt{2} Id$. The number $\sqrt{2}$ is arbitrary and any number $\gamma > 1$ leads to the same space $\mathcal{G} = \mathcal{G}_{\sqrt{2} Id}$ and its dual space $\mathcal{G}' = \mathcal{G}'_{\sqrt{2} Id}$. For $s \in \mathbb{Z}$ we simply write $\mathcal{G}_s$ instead of $\mathcal{G}_{K,s}$. These spaces were introduced and studied in [25]. Thus, we obtain that the elements $\Phi^{(n)}$, $n \in \mathbb{N}_0$, from the generalised chaos decomposition of an element $\Phi = \sum_{n=0}^{\infty} \langle \Phi^{(n)}, :\cdot^{\otimes n}: \rangle \in \mathcal{G}'_K$ satisfy $\Phi^{(n)} \in \widehat{L^2_\mathbb{C}(\mathbb{R}^n)}$. In particular, if $\Phi^{(n)} = 0$ for all but finitely many $n \in \mathbb{N}_0$ then we directly obtain $\Phi \in L^2(\mu)$. Observe that for the choice of the operator $K = \sqrt{2} Id$ the additional assumption in Remark 4.4 is obviously satisfied. Hence, in this case we formulate Theorem 2.11 as follows:

**Corollary 5.1** *Let $\Phi \in (\mathcal{N})'$. Then for $s \in \mathbb{R}$ it holds*

$$\Phi \in \mathcal{G}_s \iff \sup_{P \in \mathbb{P}} \int_{S'_\mathbb{C}(\mathbb{R})} \left| S\Phi(2^{\frac{s}{2}} P \eta) \right|^2 \nu(d\eta) < \infty.$$

*In particular, the spaces $\mathcal{G}$ and $\mathcal{G}'$ are characterised as follows:*

(i) $\Phi \in \mathcal{G} \Longleftrightarrow \forall \lambda > 0 : \sup_{P \in \mathbb{P}} \int_{S'_\mathbb{C}(\mathbb{R})} |S\Phi(\lambda P \eta)|^2 \, \nu(d\eta) < \infty$.

(ii) $\Phi \in \mathcal{G}' \Longleftrightarrow \exists \varepsilon > 0 : \sup_{P \in \mathbb{P}} \int_{S'_\mathbb{C}(\mathbb{R})} |S\Phi(\varepsilon P \eta)|^2 \, \nu(d\eta) < \infty$.

For the rest of this subsection let us fix a Brownian motion $(B_t)_{t \geq 0}$ which rises from the Kolmogorov continuity theorem as a continuous modification of the family $\left(\langle 1_{[0,t]}, \cdot \rangle\right)_{t \geq 0} \subseteq L^2(\mu)$. We denote the natural filtration of the Brownian motion $(B_t)_{t \geq 0}$ by $(\mathcal{F}_t)_{t \geq 0}$, i.e., $\mathcal{F}_t = \sigma(B_s \mid s \in [0, t])$, $t \geq 0$.

In [4,6,11,20,24] a parabolic SPDE modelling the transport of a substance in a turbulent medium is treated via white noise analysis. There the authors search for a solution $u : \mathbb{R}_+ \times \mathbb{R}^d \times \Omega : \longrightarrow \mathbb{R}$, $(t, x, \omega) \mapsto u(t, x, \omega)$ describing the concentration of the substance, where $t$ stands for the time, $x$ for the position and $\omega$ for the random parameter which will be suppressed in the following. For sake of simplicity, we only consider here the one dimensional case $d = 1$. All calculations below generalise to the multidimensional case. The SPDE under consideration is given by

$$\frac{\partial u_{t,x}}{\partial t} = \frac{1}{2}\nu(t)\frac{\partial^2 u_{t,x}}{\partial x^2} + \frac{\partial u_{t,x}}{\partial x}\sigma(t)\dot{B}_t, \quad t > 0, x \in \mathbb{R} \tag{20}$$

$$u(0, \cdot) = \delta_0, \tag{21}$$

where $\nu$ describes the molecular viscosity of the medium and $dB_t$ denotes the Itô integral w.r.t. a Brownian motion $(B_t)_{t \geq 0}$ modelling the turbulence in the medium. The initial condition (21) is a physical idealisation that at time zero the substance is only concentrated at the point $x = 0$. Thus, we obtain an analogue of an integral kernel of the SPDE (22), as known in the field of partial differential equations. Hence, more realistic and even random initial conditions can be realised via convolution, see Remark 5.7 below. In [4,24] the stochastic integral is treated in the Stratonovich sense and existence of an $L^2$-valued solution $u_{t,x}$ is shown. The Itô case is also treated in [4,24]. In [4] the solution is constructed as a generalised Brownian functional, see the last mentioned reference as well as [14] for the precise meaning. In [6,24] the solution $u_{t,x}$ is constructed in the space of Hida distributions $(\mathcal{N})'$. Furthermore, in the last mentioned reference explicit conditions on $\nu$ and $\sigma$ in terms of Hölder regularity are given, such that $u_{t,x} \in L^2(\mu)$.

In the following we use Corollary 5.1 to improve the results in [24] by giving explicit conditions on the coefficients $\nu(t)$, $\sigma(t)$ s.t. $u_{t,x} \in \mathcal{G}_s \subseteq (\mathcal{N})'$, $s \in \mathbb{Z}$, see Theorem 5.2 below.

To formulate (20), (21) in terms of white noise analysis we introduce the white noise process $(w_t)_{t \geq 0} \subseteq (S)'$. The element $w_t$ is given by its generalised chaos decomposition $w_t = \langle \delta_t, \cdot \rangle$, where $\delta_t \in S'(\mathbb{R})$ denotes the Dirac delta distribution at $t \geq 0$. Furthermore, we introduce the Wick product on $(S)'$. For $\Phi, \Psi \in (S)'$ we define the Wick product $\Phi \diamond \Psi \in (S)'$ via the $S$-transform, i.e., $S(\Phi \diamond \Psi) = S(\Phi)S(\Psi)$. Observe that the product of two $U$-functionals is again a $U$-functional, hence, $\Phi \diamond \Psi$ is well-defined by Theorem 2.8. A rigorous interpretation of (20), (21) in terms of white noise analysis is now given as follows. We search for a map $u : \mathbb{R}_+ \times \mathbb{R} \longrightarrow (\mathcal{N})'$

fulfilling

$$\frac{\partial u_{t,x}}{\partial t} = \frac{1}{2}\nu(t)\frac{\partial^2 u_{t,x}}{\partial x^2} + \sigma(t)w_t \diamond \frac{\partial u_{t,x}}{\partial x}, t > 0, x \in \mathbb{R}, \tag{22}$$

$$\left(Su_{t,\cdot}(\varphi)\right)_{t>0} \text{ is a Dirac sequence for all } \varphi \in S(\mathbb{R}). \tag{23}$$

The initial condition (23) means that for all $\varphi, g \in S(\mathbb{R})$ the following is valid

$$\lim_{t \to 0} \int_{\mathbb{R}} Su_{t,x}(\varphi)g(x)\,dx = g(0).$$

We explain the connection between the Itô term in (20) and the so called Hitsuda–Skorokhod term $\sigma(t)\frac{\partial u_{t,x}}{\partial x} \diamond w_t$ in (22) in Remark 5.7 below. In the following we denote by $\overline{\mathbb{R}}$ the extended real numbers. We formulate our existence result in the next theorem:

**Theorem 5.2** *Assume that $\nu : [0, \infty) \longrightarrow \overline{\mathbb{R}}$ is strictly positive and locally integrable and $\sigma : [0, \infty) \longrightarrow \overline{\mathbb{R}}$ is locally square integrable. If the function $(0, \infty) \ni t \mapsto \kappa(t) := \frac{\int_{[0,t]} \sigma^2(s)\,ds}{\int_{[0,t]} \nu(s)\,ds} \in \mathbb{R}$ is bounded in the vicinity of 0 then for every $T \in \mathbb{N}$ there exists an $s \in \mathbb{R}$ and a map*

$$u : (0, T] \times \mathbb{R} \longrightarrow \mathcal{G}_s$$

*satisfying (23). Furthermore, for dt-a.e. $t \in (0, T]$ and all $x \in \mathbb{R}$ the map $u$ is once differentiable w.r.t. $t$ and twice differentiable w.r.t. $x$ at $(t, x)$ and satisfies (22). In particular, for $s \in \mathbb{R}$ and $t \in (0, T]$ satisfying $2^s \kappa(t) < 1$ it holds $u_{t,x} \in \mathcal{G}_s$ for all $x \in \mathbb{R}$.*

**Proof** The same computations as in [23, Sect. 5] yield a candidate for the $S$-transform of $u$:

$$Su_{t,x}(\varphi) = \frac{1}{\sqrt{2\pi\,\vartheta(t)}} \exp\left(-\frac{1}{2\vartheta(t)}\left(x - \langle 1_{[0,t]}\sigma, \varphi\rangle\right)^2\right), \quad \varphi \in S(\mathbb{R}), \tag{24}$$

where $\vartheta(t) = \int_{[0,t]} \nu(s)\,ds, t > 0$. One easily sees that (24) defines a $U$-functional satisfying (23). Via Theorem 2.8 we obtain the element $u_{t,x} \in (\mathcal{N})'$ having $S$-transform given by (24). If $\int_{[0,t]} \nu(s)\,ds = \int_{[0,t]} \sigma(s)^2\,ds$ the corresponding Hida distribution $u_{t,x}$ is given by Donskers delta $\delta_x(\langle 1_{[0,t]}\sigma, \cdot\rangle)$, see e.g. [19, Example 13.9.]. The fact that $u_{t,x}$ satisfies (22) in the weak sense is proven [23]. In the remaining part of the proof we show that for all $T \in \mathbb{N}$ there exists an $s \in \mathbb{R}$ s.t. $u_{t,x} \in \mathcal{G}_s$ for all $(t, x) \in (0, T] \times \mathbb{R}$.

We divide the proof into two separate parts. In the first part we show that $u$ is differentiable and satisfies (22) in the above mentioned sense. In the second part we show that for all $T \in \mathbb{N}$ there exists an $s \in \mathbb{R}$ s.t. $u_{t,x} \in \mathcal{G}_s$ for all $(t, x) \in (0, T] \times \mathbb{R}$.

*Part 1* To show that $u$ is differentiable and satisfies (22) in the above mentioned sense we use [15, Theorem 4.41.]. We only show that $u$ is differentiable w.r.t. $t$ at every $(t, x) \in D \times \mathbb{R}$, where $(0, T] \backslash D$ is of Lebesgue measure zero. The treatment of the derivatives w.r.t. $x$ is easier and can be done by the same procedure. The fact that $u_{t,x}$ satisfies (22) can be seen by considering the corresponding equation for the $S$-transform $Su_{t,x}$. We make the following observation. Let $\varphi \in S(\mathbb{R})$ and $T \in \mathbb{N}$. Via the fundamental theorem of Lebesgue calculus, the functions $\vartheta$ and $t \mapsto \varrho(t) := \int_{[0,t]} \sigma(s)\varphi(s)\,ds$ are absolutely continuous and differentiable at $dt$-a.e. $t \in (0, T]$ with respective derivatives $\nu(t)$ and $\sigma(t)\varphi(t)$. The set of all $t \in (0, T]$ s.t. $\vartheta$ and $\varrho$ are differentiable at $t$ we denote by $D_1$. Hence, $Su(\cdot, x)(\varphi)$ is differentiable at $t \in D_1$ and for a zero sequence $(h_n)_{n \in \mathbb{N}}$, s.t. $|h_n| \leq \frac{t}{2}$, it holds

$$
\begin{aligned}
\frac{\partial Su_{t,x}(\varphi)}{\partial t} &= \lim_{n \to \infty} \frac{Su(t + h_n, x)(\varphi) - S(u_{t,x}(\varphi)}{h_n} \\
&= \frac{1}{2\vartheta(t)^{\frac{3}{2}}\sqrt{2\pi}} exp\left(-\frac{1}{2\vartheta(t)}\left(x - \int_{[0,t]}\sigma(s)\varphi(s)\,ds\right)^2\right) \\
&\quad \times \left(-\nu(t) - 2\sigma(t)\varphi(t)\left(x - \int_{[0,t]}\sigma(s)\varphi(s)\,ds\right)\right. \\
&\quad \left. -\frac{\nu(t)\left(x - \int_{[0,t]}\sigma(s)\varphi(s)\,ds\right)^2}{\vartheta(t)}\right).
\end{aligned}
$$

Hence, for $z \in \mathbb{C}$ we obtain the estimate

$$
\begin{aligned}
\left|\frac{\partial Su_{t,x}(z\varphi)}{\partial t}\right| &\leq \frac{1}{2\vartheta(t)^{\frac{3}{2}}\sqrt{2\pi}}\exp\left(\frac{\|\sigma\|_{L_T^2}^2}{\vartheta(t)}|z|^2\|\varphi\|_{L_T^2}^2\right) \\
&\quad \times \left(|\nu(t)| + |z|^2\|\varphi\|_\infty^2\,\sigma(t)^2 + \left(|x| + \|\sigma\|_{L_T^2}|z|\|\varphi\|_{L_T^2}\right)^2\right. \\
&\quad \left. +\frac{|\nu(t)|\left(|x| + \|\sigma\|_{L_T^2}|z|\|\varphi\|_{L_T^2}\right)^2}{\vartheta(t)}\right) \\
&\leq \frac{1}{2\vartheta(t)^{\frac{3}{2}}\sqrt{2\pi}}\exp\left(\frac{\|\sigma\|_{L_T^2}^2}{\vartheta(t)}|z|^2\|\varphi\|_{L_T^2}^2\right) \\
&\quad \times \left(|\nu(t)| + \exp(|z|^2\|\varphi\|_\infty^2)\sigma(t)^2 + C_1\exp(\|\sigma\|_{L_T^2}^2|z|^2\|\varphi\|_{L_T^2}^2)\right. \\
&\quad \left. +|\nu(t)|\frac{C_1\exp\left(\|\sigma\|_{L_T^2}^2|z|^2\|\varphi\|_{L_T^2}^2\right)}{\vartheta(t)}\right) \\
&\leq C_3(t)\exp\left(C_2(t)|z|^2\|\varphi\|_p^2\right)\left(|\nu(t)| + \sigma^2(t) + 1\right)
\end{aligned}
$$

where $\|\cdot\|_{L_T^2}$ denotes the $L^2((0, T))$-norm, $\|\cdot\|_\infty$ the $L^\infty(\mathbb{R})$-norm, $p \in \mathbb{N}$ is chosen s.t. for all $\varphi \in S(\mathbb{R})$ it holds $\max\left\{\|\varphi\|_{L_T^2}, \|\varphi\|_\infty\right\} \leq \|\varphi\|_p := \|A^p\varphi\|_{L^2(\mathbb{R})}$ and

$$C_1 = \max\{2|x|, 2\}, \quad C_2(t) = \frac{1}{\vartheta(t)}\|\sigma\|_{L_T^2}^2 + 1 + \|\sigma\|_{L_T^2}^2,$$

$$C_3(t) = \frac{\max\left\{C_1, 1 + \frac{C_1}{\vartheta(t)}\right\}}{2\vartheta(t)^{\frac{3}{2}}\sqrt{2\pi}}.$$

Observe that $C_2$ and $C_3$ are decreasing. Applying the fundamental theorem of Lebesgue calculus to $Su(\cdot, x)(\varphi)$ it holds

$$\left|\frac{Su(t + h_n, x)(\varphi) - S(u_{t,x}(\varphi))}{h_n}\right| = \left|\frac{1}{h_n}\int_{[t, t+h_n]} \frac{\partial Su(s, x)(\varphi)}{\partial s}\, ds\right|$$

$$\leq C_3\left(\frac{t}{2}\right)\exp\left(C_2\left(\frac{t}{2}\right)|z|^2\|\varphi\|_p^2\right)\frac{1}{|h_n|}\int_{[t, t+h_n]} |v(s)| + \sigma^2(s) + 1\, ds$$

Via the fundamental theorem of Lebesgue calculus it holds for $dt$-a.e. $t$

$$C_4(t) := \sup_{n \in \mathbb{N}} \frac{1}{|h_n|}\int_{[t, t+h_n]} |v(s)| + \sigma^2(s) + 1\, ds < \infty. \tag{25}$$

We denote the set of all $t \in (0, T]$ s.t. (25) holds true by $D_2$. We conclude that for $t \in D := D_1 \cap D_2$ it holds

$$\lim_{n \to \infty} \frac{Su(t + h_n, x)(\varphi) - S(u_{t,x}(\varphi)}{h_n} \text{ exists,}$$

$$\left|\frac{Su(t + h_n, x)(\varphi) - S(u_{t,x}(\varphi)}{h_n}\right| \leq C_3\left(\frac{t}{2}\right)C_4(t)\exp\left(C_2\left(\frac{t}{2}\right)|z|^2\|\varphi\|_p^2\right).$$

Now we apply [15, Theorem 4.41.] and obtain that $u$ is differentiable w.r.t. $t$ at $(t, x) \in D \times \mathbb{R}$. The first part is finished.

*Part 2* It is left to show that for every time $T \in \mathbb{N}$ we can find an $s \in \mathbb{R}$ s.t. $u_{t,x} \in \mathcal{G}_s$ for all $t \in (0, T]$ and $x \in \mathbb{R}$. We prove this by using Corollary 5.1. To this end let $P \in \mathbb{P}$ be a projection as in Definition 2.10, $\eta_1 + i\eta_2 = \eta \in S'_{\mathbb{C}}(\mathbb{R})$, where $\eta_1, \eta_2 \in S'(\mathbb{R})$ and $\varepsilon > 0$.

We obtain

$$|Su_{t,x}(\varepsilon P\eta)|^2 = \left|\frac{1}{\sqrt{2\pi\vartheta(t)}}\exp\left(-\frac{1}{2\vartheta(t)}\left(x - \varepsilon\left\langle 1_{[0,t]}\sigma, P\eta\right\rangle\right)^2\right)\right|^2$$

$$= \left|\frac{1}{\sqrt{2\pi\vartheta(t)}}\exp\left(-\frac{1}{2\vartheta(t)}\left(x - \varepsilon\left\langle P(1_{[0,t]}\sigma), \eta\right\rangle\right)^2\right)\right|^2$$

$$= \frac{1}{2\pi\vartheta(t)}\exp\left(-\frac{\varepsilon^2}{\vartheta(t)}\left(\frac{x}{\varepsilon} - \left\langle P(1_{[0,t]}\sigma), \eta_1\right\rangle\right)^2\right)$$

$$\exp\left(\frac{\varepsilon^2}{\vartheta(t)}\left\langle P(1_{[0,t]}\sigma),\eta_2\right\rangle^2\right).$$

Now we calculate the integral of $\left|S\left(u_{t,x}\right)(\varepsilon P\cdot)\right|^2$ w.r.t. the measure $\nu=\mu_{\frac12}\otimes\mu_{\frac12}$. Observe that the law $\mu_{\frac12}\circ\left\langle\frac{P(1_{[0,t]}\sigma)}{\|P(1_{[0,t]}\sigma)\|_{L^2(\mathbb{R})}},\cdot\right\rangle^{-1}$ is the centered Gaussian measure with variance $\frac12$. Thus we conclude

$$
\begin{aligned}
&\int_{S'_{\mathbb{C}}(\mathbb{R})}\left|S\left(u_{t,x}\right)(\varepsilon P\eta)\right|^2\,d\nu(\eta)\\
&=\frac{1}{\sqrt{2\pi\,\vartheta(t)\pi}}\int_{\mathbb{R}}\exp\left(-\frac{\varepsilon^2}{\vartheta(t)}\left(\frac{x}{\varepsilon}-y_1\left\|P(1_{[0,t]}\sigma)\right\|_{L^2(\mathbb{R})}\right)^2\right)\exp\left(-|y_1|^2\right)\,dy_1\\
&\quad\times\frac{1}{\sqrt{2\pi\,\vartheta(t)\pi}}\int_{\mathbb{R}}\exp\left(\frac{\varepsilon^2}{\vartheta(t)}y_2^2\left\|P(1_{[0,t]}\sigma)\right\|_{L^2(\mathbb{R})}^2\right)\exp\left(-|y_2|^2\right)\,dy_2\\
&=\frac{1}{\sqrt{2\pi\,\vartheta(t)\pi}}\int_{\mathbb{R}}\exp\left(-\left(1+\frac{\varepsilon^2\left\|P(1_{[0,t]}\sigma)\right\|_{L^2(\mathbb{R})}^2}{\vartheta(t)}\right)y_1^2\right.\\
&\quad\left.+y_1\frac{2\varepsilon x\left\|P(1_{[0,t]}\sigma)\right\|_{L^2(\mathbb{R})}}{\vartheta(t)}-\frac{x^2}{\vartheta(t)}\right)\,dy_1\\
&\quad\times\frac{1}{\sqrt{2\pi\,\vartheta(t)\pi}}\int_{\mathbb{R}}\exp\left(-\left(1-\frac{\varepsilon^2\left\|P(1_{[0,t]}\sigma)\right\|_{L^2(\mathbb{R})}^2}{\vartheta(t)}\right)y_2^2\right)\,dy_2. \quad (26)
\end{aligned}
$$

From this point we see that a necessary condition for $\sup_{P\in\mathbb{P}}\int_{S_{\mathbb{C}}(\mathbb{R})}\left|S\left(u_{t,x}\right)(\varepsilon P\eta)\right|^2$ $\nu(d\eta)$ to be finite for some $\varepsilon>0$ is that $1\pm\frac{\varepsilon^2\|P(1_{[0,t]}\sigma)\|_{L^2(\mathbb{R})}^2}{\vartheta(t)}>0$. Since $\frac{\varepsilon^2\|P(1_{[0,t]}\sigma)\|_{L^2(\mathbb{R})}^2}{\vartheta(t)}\le\varepsilon^2\kappa(t)$ we choose $\varepsilon>0$ s.t. $0<\varepsilon^2\kappa(t)<1$.

Now we can evaluate the Gaussian integrals in (26) and we seek

$$
\begin{aligned}
&\int_{S'_{\mathbb{C}}(\mathbb{R})}\left|S\left(u_{t,x}\right)(\varepsilon P\eta)\right|^2\,d\nu(\eta)\\
&=\frac{\exp\left(-\frac{x^2}{\vartheta(t)}\right)}{2\pi\,\vartheta(t)}\frac{1}{\sqrt{1-\frac{\varepsilon^4\|P(1_{[0,t]}\sigma)\|_{L^2(\mathbb{R})}^4}{\vartheta(t)^2}}}\exp\left(\frac{\varepsilon^2\left\|P(1_{[0,t]}\sigma)\right\|_{L^2(\mathbb{R})}^2x^2}{\vartheta(t)\varepsilon^2\left\|P(1_{[0,t]}\sigma)\right\|_{L^2(\mathbb{R})}^2+\vartheta(t)^2}\right).
\end{aligned}
$$

We conclude by Corollary 4.3 that for $0<\varepsilon^2\kappa(t)<1$ it holds

$$\sup_{P\in\mathbb{P}}\int_{S'_{\mathbb{C}}(\mathbb{R})}\left|S\left(u_{t,x}\right)(\varepsilon P\eta)\right|^2\nu(d\eta)$$

$$= \frac{\exp\left(-\frac{x^2}{\vartheta(t)}\right)}{2\pi\,\vartheta(t)} \frac{1}{\sqrt{1 - \frac{\varepsilon^4 \|1_{[0,t]}\sigma\|_{L^2(\mathbb{R})}^4}{\vartheta(t)^2}}} \exp\left(\frac{\varepsilon^2 \|1_{[0,t]}\sigma\|_{L^2(\mathbb{R})}^2 x^2}{\vartheta(t)\varepsilon^2 \|1_{[0,t]}\sigma\|_{L^2(\mathbb{R})}^2 + \vartheta(t)^2}\right).$$

Observe that by the assumptions on the coefficients $\nu$ and $\sigma$ it holds that $(0, \infty) \ni t \mapsto \kappa(t) \in \mathbb{R}$ is continuous. Consequently, by assumption $\kappa$ is bounded on finite intervals $(0, T]$, $T \in \mathbb{N}$. Hence, for every $T \in \mathbb{N}$ we can choose $\varepsilon > 0$ s.t. $0 < \varepsilon^2 \kappa(t) < 1$ for all $t \in (0, T]$. Eventually, we conclude by (20) that for $s \in \mathbb{R}$ satisfying $2^{\frac{s}{2}} = \varepsilon$ it holds

$$u_{t,x} \in \mathcal{G}_s, \text{ for all } t \in (0, T], x \in \mathbb{R}.$$

In the case $\kappa(t) < 1$ we can choose $\varepsilon \geq 1$ which implies $u_{t,x} \in L^2(\mu)$. $\qquad\square$

**Remark 5.3** (i) As observed in [24], the solution $u_{t,x}$ passes through Donskers delta $\delta_x(\langle 1_{[0,t]}\sigma, \cdot\rangle)$ each time $\kappa(t)$ passes through the value 1. In particular, the solution $u$ satisfies

$$u_{t,x} \begin{cases} \in L^2(\mu), & \text{if } \kappa(t) < 1, \\ = \delta_x(\langle 1_{[0,t]}\sigma, \cdot\rangle), & \text{if } \kappa(t) = 1, \\ \in \mathcal{G}', & \text{if } \kappa(t) > 1. \end{cases}$$

In the case $\kappa(t) < 1$ the solution $u_{t,x}$ is explicitly given by, see [24, Eq. (3.13)],

$$u_{t,x} = \left(2\pi \int_{[0,t]} \nu(s) - \sigma^2(s)\,ds\right)^{-\frac{1}{2}}$$
$$\exp\left(-\left(2\int_{[0,t]} \nu(s) - \sigma^2(s)\,ds\right)^{-1}(x - \langle 1_{[0,t]}\sigma, \cdot\rangle)^2\right).$$

In particular, in the case $\kappa(t) < 1$ Lemma 5.6 below implies that a version of $u_{t,x}$ is measurable w.r.t. $\mathcal{F}_t$.

(ii) Observe that the calculation in the proof of Theorem 5.2 shows that Donsker's delta $\delta_x(\langle f, \cdot\rangle)$, $x \in \mathbb{R}$, $f \in L^2(\mathbb{R})$, is an element of $\mathcal{G}_s$ for all $s \in (-\infty, 0)$ and $\delta_x(\langle f, \cdot\rangle) \notin L^2(\mu)$.

In the following let $\Phi = \sum_{n=0}^{\infty} \langle \Phi^{(n)}, ::^{\otimes n} :\rangle \in \mathcal{G}'$. Recall that the $S$-transform of $\Phi$ is given by

$$S\Phi(\varphi) = \sum_{n=0}^{\infty} \langle \varphi^{\otimes n}, \Phi^{(n)} \rangle \qquad (27)$$

and the dual pairing $\langle \varphi^{\otimes n}, \Phi^{(n)} \rangle$ is given by the scalar product $\left(\varphi^{\otimes n}, \overline{\Phi^{(n)}}\right)_{L^2(\mathbb{R}^n)}$. Hence, the $S$-transform of $\Phi$ admits a natural extension to $L^2(\mathbb{R})$ which is also given by (27). The next lemma follows immediately from the polarisation identity.

**Lemma 5.4** *Let $\Phi = \sum_{n=0}^{\infty} \langle \Phi^{(n)}, : \cdot^{\otimes n} : \rangle \in \mathcal{G}'$ and $I \subseteq \mathbb{R}$ be measurable. Then the following are equivalent:*

(i) $supp(\Phi^n) \subseteq I^n$ *for all $n \in \mathbb{N}$.*
(ii) $S\Phi(\varphi) = S\Phi(1_I \varphi)$ *for all $\varphi \in S(\mathbb{R})$.*

**Lemma 5.5** *Let $F = \sum_{n=0}^{\infty} \langle F^{(n)}, : \cdot^{\otimes n} : \rangle \in L^2(\mu)$ and $I \subseteq \mathbb{R}$. Then there exists a version $\tilde{F}$ of $F$ which is measurable w.r.t. $\sigma(\langle \xi, \cdot \rangle, \xi \in S(\mathbb{R}), supp(\xi) \subseteq I)$ if and only if $supp(F^{(n)}) \subseteq I^n$ for all $n \in \mathbb{N}$.*

**Proof** It suffices to prove the statement for $F = \langle F^{(n)}, : \cdot^{\otimes n} : \rangle$, $n \in \mathbb{N}$. Necessity can be proven as in [14, Proposition 4.5.]. Sufficiency follows from construction of the element $\langle F^{(n)}, : \cdot^{\otimes n} : \rangle \in L^2(\mu)$. $\qquad\square$

For the next lemma recall the Brownian motion $(B_t)_{t \geq 0}$ as well as the corresponding natural filtration $(\mathcal{F}_t)_{t \geq 0}$ defined above.

**Lemma 5.6** *Assume that $F : S'(\mathbb{R}) \longrightarrow \mathbb{C}$ is measurable w.r.t. $\mathcal{F}_t$, $t \in [0, \infty)$. Then there exists a $G : S'(\mathbb{R}) \longrightarrow \mathbb{C}$ which is measurable w.r.t. $\mathcal{A}_t := \sigma(\langle \xi, \cdot \rangle \mid \xi \in S(\mathbb{R}), supp(\xi) \subseteq [0, t))$ and it holds $F = G$ $\mu$-a.e. and vice versa.*

**Proof** Using [3, Corollary 2.9] it suffices to show that a version of $B_s$ is measurable w.r.t. $\mathcal{A}_t$ for $s \in [0, t]$ and that a version of $\langle \xi, \cdot \rangle$ is measurable w.r.t. $\mathcal{F}_t$ for $\xi \in S(\mathbb{R})$ with $\text{supp}(\xi) \subseteq [0, t)$. Both statements follow by the construction of $\langle 1_{[0,s]}, \cdot \rangle \in L^2(\mu)$, $s \in [0, t]$. $\qquad\square$

Several remarks are in order.

**Remark 5.7** (i) Let $f : \mathbb{R} \longrightarrow \mathcal{G}'$ be given s.t. for $(t, x) \in (0, \infty) \times \mathbb{R}$ the map

$$\mathbb{R} \ni y \mapsto f(y) \diamond u(t, x - y) \in \mathcal{G}'$$

is weakly in $L^1(\mathbb{R}, \mathcal{B}, dx)$, where $\mathcal{B}$ is the Borel $\sigma$-field and $dx$ the Lebesgue measure on $\mathbb{R}$, respectively and $u$ is defined in (24). I.e., for every $F \in \mathcal{G}$ holds that $\mathbb{R} \ni y \mapsto \langle F, f(y) \diamond u(t, x - y) \rangle \in \mathbb{C}$ is in $L^1(\mathbb{R}, \mathcal{B}, dx)$. Then we can define the Pettis-integral

$$u_{t,x}^f := \int_{\mathbb{R}} f(y) \diamond u_{t,x-y} \, dy \in \mathcal{G}', \tag{28}$$

see also [15, Proposition 8.1] and [25, Proposition 2.6]. Under additional assumptions on $f$ we obtain that $u^f$ satisfies the initial condition $u_{0,x}^f = \lim_{t \to 0} u_{t,x}^f = f(x)$ for all $x \in \mathbb{R}$. If we assume that the time and space derivatives $\frac{d}{dt}$, $\frac{d}{dx}$ and $\frac{d^2}{dx^2}$ commute with the Pettis-integral (28) we obtain that $u^f : \mathbb{R}_{\geq 0} \times \mathbb{R} \longrightarrow \mathcal{G}'$ satisfies (22) with the initial condition $u_{0,x}^f = f(x)$.

(ii) If we assume for the sake of simplicity that the initial data $f$ in (i) is deterministic and an element from $S(\mathbb{R})$, all steps in (i) are justified and we obtain that $u^f$ given by (28) is a solution to (22) with initial condition $u_{(}^f 0, x) = f(x)$ for all $x \in \mathbb{R}$.

From Theorem 5.2 and Lemma 5.4 we can conclude that $(u^f_{t,x})_{t \in [0,T]}$, $x \in \mathbb{R}$, is a generalised stochastic process and adapted in the sense of [1, Definition 1]. In particular, if $u^f_{t,x} \in L^2(\mu)$ then it holds by Lemmas 5.5 and 5.6 that a version of $u^f_{t,x}$ is measurable w.r.t. $\mathcal{F}_t$ for all $x \in \mathbb{R}$. It is well-known that for such a process the Itô integral and the Hitsuda–Skorokhod integral coincide, see e.g. [15, Theorem 8.7.]. Indeed, if for $x \in \mathbb{R}$ and $T \in [0,\infty)$ the process $\left( \sigma(s) \frac{\partial u^f_{s,x}}{\partial x} \right)_{s \in [0,T]}$ is in $L^2([0,T]; L^2(\mu))$ then for $t \in [0,T]$ the following identity is valid

$$\int_0^t \sigma(s) \frac{\partial u^f_{s,x}}{\partial x} \diamond w_s \, ds = \int_0^t \sigma(s) \frac{\partial u^f_{s,x}}{\partial x} \, dB_s \quad \mu\text{-a.e.,}$$

## 5.2 Stochastic Heat equation with general multiplicative colored Noise

In this subsection we apply our characterisation theorem to derive new results regarding the stochastic heat equation with multiplicative noise having a very general covariance structure. The results in this subsection are an extension of the results in [17, Sect. 3,4]. Indeed, we combine our main result with the calculations made in [17] to improve the results given in the last mentioned reference. We start by recalling the central objects introduced in Sects. 3 and 4 of [17] and relate them to our general functional analytic framework. In Theorems 5.10 and 5.12 we formulate our new results regarding the heat equation with multiplicative stochastic source term in the case of Skorokhod and Stratonovich product, respectively.

Before we begin with our considerations let us introduce some notation. Throughout this entire subsection $d \in \mathbb{N}$ is fixed and by $p_t$, $t > 0$, we denote the heat kernel given by $p_t(x) = (2\pi t)^{-\frac{d}{2}} \exp(-\frac{|x|^2}{2t})$, $x \in \mathbb{R}^d$. Furthermore, let $(\Omega, \mathcal{A}, \boldsymbol{P})$ be an arbitrary probability space carrying two independent $d$-dimensional Brownian motions $(B_t)_{t \geq 0}$, $(\tilde{B}_t)_{t \geq 0}$. For $x \in \mathbb{R}^d$ we denote by $(B^x_t)_{t \geq 0}$ and $(\tilde{B}^x_t)_{t \geq 0}$ the processes $(B_t + x)_{t \geq 0}$ and $(\tilde{B}_t + x)_{t \geq 0}$, respectively.

The stochastic partial differential equation we consider is informally given by

$$\frac{\partial u_{t,x}}{\partial t} = \frac{1}{2} \Delta u_{t,x} + u_{t,x} \dot{W}_{t,x}, \quad t > 0, x \in \mathbb{R}^d, \tag{29}$$

$$u_{0,x} = u_0(x), \quad x \in \mathbb{R}^d, \tag{30}$$

with continuous and bounded initial condition $u_0$. The product between $u_{t,x}$ and the centered Gaussian process $\dot{W}_{t,x}$, $t > 0$, $x \in \mathbb{R}^d$, is treated in the Skorokhod and the Stratonovich case, see [17]. The covariance structure of $\dot{W}$ is given by

$$\mathbb{E}\left[ \dot{W}_{t,x} \dot{W}_{s,y} \right] = \gamma(t - s) \Lambda(x - y), \tag{31}$$

where $\gamma$ and $\Lambda$ are generalized functions. The rigorous interpretation of (29) and (31) is given below. We work under the following assumptions on $\gamma$ and $\Lambda$.

**Assumption 5.8** Let $\gamma : \mathbb{R} \longrightarrow \mathbb{R}_+$, $\Lambda : \mathbb{R}^d \longrightarrow \mathbb{R}_+$ be measurable, non-negative definite functions, s.t., $\gamma \in S'(\mathbb{R}) \cap L^1_{loc}(\mathbb{R})$ and $\Lambda \in S'(\mathbb{R}^d)$, where $S'(\mathbb{R})$ and $S'(\mathbb{R}^d)$ are the spaces of tempered distributions over $\mathbb{R}$ and $\mathbb{R}^d$, respectively. The Fourier transforms $\rho = \mathcal{F}\gamma$ and $\sigma = \mathcal{F}\Lambda$ are tempered measures on $\mathbb{R}$ and $\mathbb{R}^d$, respectively and the product measure $\rho \otimes \sigma$ has full topological support. The measure $\sigma$ satisfies

$$\int_{\mathbb{R}^d} \frac{1}{1 + |\xi|^2} \sigma(d\xi) < \infty.$$

To specify our functional analytic framework let $\gamma$ and $\Lambda$ satisfy Assumption 5.8. Throughout this subsection we choose $\mathcal{N} = S(\mathbb{R}^{d+1})$ the real-valued Schwartz functions and define $\mathcal{H}$ as the abstract completion of $\mathcal{N}$ w.r.t. the inner product

$$(f, g)_{\mathcal{H}} = \int_{\mathbb{R}^2} \int_{\mathbb{R}^{2d}} f(t, x) \overline{g(s, y)} \gamma(t - s) \Lambda(x - y) \, dx \, dy \, dt \, ds$$
$$= \int_{\mathbb{R}} \int_{\mathbb{R}^d} \mathcal{F}f(\alpha, \xi) \overline{\mathcal{F}g(\alpha, \xi)} \rho(d\alpha) \sigma(d\xi).$$

In particular, $\mu$ denotes the mean zero Gaussian measure defined on $\mathcal{N}' = S'(\mathbb{R}^{d+1})$ with covariance given by the inner product $(\cdot, \cdot)_{\mathcal{H}}$. For more details see Example 2.12(ii). As in the previous section we consider $K = \sqrt{2}Id$ and for $s \in \mathbb{R}$ we denote the spaces $\mathcal{G}_K, \mathcal{G}_{K,s}$ simply by $\mathcal{G}, \mathcal{G}_s$, respectively. Equivalently as in Corollary 5.1 our main result can be formulated under the above mentioned framework as

$$\Phi \in \mathcal{G}_s \iff \sup_{P \in \mathbb{P}} \int_{S'_\mathbb{C}(\mathbb{R}^{d+1})} \left| S\Phi(2^{\frac{s}{2}} P\eta) \right|^2 \nu(d\eta) < \infty,$$

where $s \in \mathbb{R}$. Before we proceed, some comments on Assumption 5.8 in comparison to the Assumptions in the reference [17] are necessary.

*Remark 5.9* (i) The assumptions made in Assumption 5.8 on $\gamma$ and $\Lambda$ are stronger than the ones made in [17, Theorem 3.6.]. In particular, we assume that $\gamma$ is a tempered distribution, which is a stronger assumption than assuming merely local integrability of $\gamma$. The assumption on the support of $\rho \otimes \sigma = \mathcal{F}\gamma \otimes \mathcal{F}\Lambda$ guarantees that the bilinear form $(\cdot, \cdot)_{\mathcal{H}}$ is positive definite on $\mathcal{N} = S(\mathbb{R}^{d+1})$. Furthermore, this property of $\gamma$ and $\Lambda$ implies that the nuclear space $\mathcal{N}$ is continuously embedded into $\mathcal{H}$, which is necessary to apply our results from the previous sections.

(ii) One can also drop the assumption that $\rho \otimes \sigma$ have full topological support in Assumption 5.8 by considering instead of $\mathcal{N} = S(\mathbb{R}^{d+1})$ the corresponding nuclear quotient space $\mathcal{N} = S(\mathbb{R}^{d+1}) \bmod N_0$, where $N_0 = \{f \in S(\mathbb{R}^{d+1}) \mid \int_{\mathbb{R}^{d+1}} |\mathcal{F}f|^2 \, d\rho \otimes \sigma = 0\}$. To simplify our considerations below we stick to the assumptions as stated in Assumption 5.8.

(iii) Important examples of $\gamma$ and $\Lambda$ can be found in [16,17], corresponding for example to white or fractional noise.

**5.2.1 Skorokhod case**

To treat the Skorokhod case let us recall what has been achieved in [17]. The authors in [17] investigate the Skorokhod case of (29), (30) by considering an approximate equation given by

$$\frac{\partial u_{t,x}^{\varepsilon,\delta}}{\partial t} = \frac{1}{2}\Delta u_{t,x}^{\varepsilon,\delta} + u_{t,x}^{\varepsilon,\delta} \diamond \dot{W}_{t,x}^{\varepsilon,\delta}, \quad t > 0, x \in \mathbb{R}^d, \tag{32}$$

$$u_{0,x}^{\varepsilon,\delta} = u_0(x), \quad x \in \mathbb{R}^d, \tag{33}$$

where $\diamond$ denotes the Wick product, see [19, Definition 8.11.] as well as [17, Equation (2.12)], and $\dot{W}_{t,x}^{\varepsilon,\delta}$ is given by the element with only first order chaos decomposition $\dot{W}_{t,x}^{\varepsilon,\delta} = \left\langle w_{t,x}^{\varepsilon,\delta}, \cdot \right\rangle$, where $w_{t,x}^{\varepsilon,\delta} \in \mathcal{H}$ is given by

$$w_{t,x}^{\varepsilon,\delta}(s, y) = \frac{1}{\delta} 1_{[0,\delta]}(t - s) 1_{[0,t]}(s) p_\varepsilon(x - y), \quad s > 0, y \in \mathbb{R}^d.$$

By [17, Eq. (3.18)] the mild solution to the approximate equation (32), (33) is given by the Bochner integral in $L^p(\mathcal{N}', \mu)$, $p \in [1, \infty)$,

$$u_{t,x}^{\varepsilon,\delta} = \mathbb{E}\left[u_0(B_t^x) : \exp\left(\left\langle A_{t,B^x}^{\varepsilon,\delta}, \cdot \right\rangle\right) : \right],$$

where $:\exp\left(\left\langle A_{t,B^x}^{\varepsilon,\delta}, \cdot \right\rangle\right):$ denotes the Wick exponential of $A_{t,B^x}^{\varepsilon,\delta} \in \mathcal{H}$ which is given by

$$A_{t,B^x}^{\varepsilon,\delta}(r, y) = \frac{1}{\delta}\left(\int_0^{\delta \wedge (t-r)} p_\varepsilon(B_{t-r-s}^x - y)\, ds\right) 1_{[0,t]}(r), \quad r \in \mathbb{R}, y \in \mathbb{R}^d,$$

and $\mathbb{E}$ denotes integration w.r.t. $\mathbb{P}$ in the sense of Bochner. Consequently, the S-transform of $u_{t,x}^{\varepsilon,\delta}$ at $\varphi \in \mathcal{N}_{\mathbb{C}}$ is given by

$$Su_{t,x}^{\varepsilon,\delta}(\varphi) = \mathbb{E}\left[u_0(B_t^x) S\left(:\exp\left(\left\langle A_{t,x}^{\varepsilon,\delta}, \cdot \right\rangle\right):\right)(\varphi)\right]$$

$$= \mathbb{E}\left[u_0(B_t^x) \exp\left((A_{t,x}^{\varepsilon,\delta}, \varphi)_{\mathcal{H}}\right)\right].$$

In [17, Theorem 3.6.] it is shown that the limit $u_{t,x} = \lim_{\varepsilon \to 0}\lim_{\delta \to 0} u_{t,x}^{\varepsilon,\delta}$ exists for all $t > 0$, $x \in \mathbb{R}^d$ in $L^p(\mathcal{N}', \mu)$ for all $p \in [1, \infty)$ and $(u_{t,x})_{t>0,x\in\mathbb{R}^d}$ coincides with the mild solution of (29), (30), see also the last mentioned reference for the definition of mild solution. With this notation in mind we can formulate our results for the Skorokhod case in the following theorem.

**Theorem 5.10** *Let $u_0 : \mathbb{R}^d \longrightarrow \mathbb{R}$ be continuous and bounded and $\gamma$ and $\Lambda$ satisfy Assumptions 5.8. For every $t > 0$, $x \in \mathbb{R}^d$ the element $u_{t,x} = \lim_{\varepsilon \to 0}\lim_{\delta \to 0} u_{t,x}^{\varepsilon,\delta}$ is* *contained in $\mathcal{G}$ and for $\lambda \in (0, \infty)$ the $\mathcal{G}_\lambda$-norm of $u_{t,x}$ can be estimated by*

$$\left\| u_{t,x} \right\|_\lambda^2 \leq \mathbb{E}\left[ u_0(B_t^x) u_0(\tilde{B}_t^x) \exp\left( 4\lambda^2 \int_0^t \int_0^t \gamma(r-s) \Lambda(B_r^x - \tilde{B}_s^x)\, ds\, dr \right) \right],$$

(34)

*where $\mathbb{E}$ denotes integration w.r.t. $\mathbb{P}$.*

**Proof** Let $t > 0$, $x \in \mathbb{R}^d$ be arbitrary and $\lambda \in (0, \infty)$. To show $u_{t,x} \in \mathcal{G}_\lambda$ it suffices to show that the norm $\left\| u_{t,x}^{\varepsilon,\delta} \right\|_\lambda$ is uniformly bounded in $\varepsilon, \delta > 0$, due to the Banach-Alaoglu theorem. To show this, we use our main result. Let us also fix $\varepsilon, \delta > 0$ and let $(e_i)_{i\in\mathbb{N}}$ denote a real orthonormal basis of $\mathcal{H}$ which is contained in $\mathcal{N}$. For $m \in \mathbb{N}$, we denote by $P_m$ the projection given by $P_m = \sum_{j=1}^m \langle e_j, \cdot \rangle e_j$. From Assumption 5.8 we conclude that the real random variable $\left\| A_{t,B^x}^{\varepsilon,\delta} \right\|_\mathcal{H}$ is bounded. Thus, we can use Fubinis Theorem and obtain

$$\int_{S'_\mathbb{C}(\mathbb{R}^{d+1})} \left| S u_{t,x}^{\varepsilon,\delta}(\lambda P_m \eta) \right|^2 \nu(d\eta)$$

$$= \mathbb{E}\left[ u_0(B_t^x) u_0(\tilde{B}_t^x) \int_{S'_\mathbb{C}(\mathbb{R}^{d+1})} \exp\left( \lambda \left\langle P_m \eta, A_{t,B^x}^{\varepsilon,\delta} \right\rangle + \lambda \overline{\left\langle P_m \eta, A_{t,\tilde{B}^x}^{\varepsilon,\delta} \right\rangle} \right) \right]$$

$$= \mathbb{E}\left[ u_0(B_t^x) u_0(\tilde{B}_t^x) \int_{S'_\mathbb{C}(\mathbb{R}^{d+1})} \exp\left( \lambda \left\langle \eta, P_m A_{t,B^x}^{\varepsilon,\delta} \right\rangle + \lambda \overline{\left\langle \eta, P_m A_{t,\tilde{B}^x}^{\varepsilon,\delta} \right\rangle} \right) \nu(d\eta) \right].$$

Proposition 3.1 yields

$$\int_{S'_\mathbb{C}(\mathbb{R}^{d+1})} \exp\left( \lambda \left\langle \eta, P_m A_{t,B^x}^{\varepsilon,\delta} \right\rangle + \lambda \overline{\left\langle \eta, P_m A_{t,\tilde{B}^x}^{\varepsilon,\delta} \right\rangle} \right) \nu(d\eta)$$

$$= \exp\left( 4\lambda^2 (P_m A_{t,B^x}^{\varepsilon,\delta}, P_m A_{t,\tilde{B}^x}^{\varepsilon,\delta})_\mathcal{H} \right).$$

Observe that the random variables $X^m := (P_m A_{t,B^x}^{\varepsilon,\delta}, P_m A_{t,\tilde{B}^x}^{\varepsilon,\delta})_\mathcal{H}$, $m \in \mathbb{N}$, are uniformly bounded. Furthermore, we have

$$X^m \xrightarrow{m\to\infty} \left( A_{t,B^x}^{\varepsilon,\delta}, A_{t,\tilde{B}^x}^{\varepsilon,\delta} \right)_\mathcal{H}$$

$$= \int_0^t \int_0^t \gamma(r-s) \frac{1}{\delta^2} \int_0^{\delta\wedge(t-r)} \int_0^{\delta\wedge(t-s)} \int_{\mathbb{R}^d} \int_{\mathbb{R}^d}$$

$$\times \Lambda(y_1 - y_2) p_\varepsilon(B_{t-r-t_1}^x - y_1) p_\varepsilon(B_{t-s-t_2}^x - y_2)\, dy_1\, dy_2\, dt_1\, dt_2\, dr\, ds$$

$$\xrightarrow{\delta\to 0} \int_0^t \int_0^t \gamma(r-s) \int_{\mathbb{R}^d} \int_{\mathbb{R}^d}$$

$$\times \Lambda(y_1 - y_2) p_\varepsilon(B_{t-r}^x - y_1) p_\varepsilon(B_{t-s}^x - y_2)\, dy_1\, dy_2\, dr\, ds$$

$$\xrightarrow{\varepsilon \to 0} \int_0^t \int_0^t \gamma(r-s) \Lambda(B_{t-s} - \tilde{B}_{t-s}) \, dr \, ds \in L^1(\Omega, \mathbb{P}),$$

where the convergence takes place in $L^1(\Omega, \mathbb{P})$, see [17, proof of Theorem 3.6]. The dominated convergence theorem implies

$$\sup_{m \in \mathbb{N}} \int_{S'_{\mathbb{C}}(\mathbb{R}^{d+1})} \left| S u_{t,x}^{\varepsilon, \delta}(\lambda P_m \eta) \right|^2 \nu(d\eta) < \infty.$$

From Corollary 4.3 we conclude

$$\left\| u_{t,x}^{\varepsilon, \delta} \right\|_{\lambda}^2 = \lim_{m \to \infty} \mathbb{E} \left[ u_0(B_t^x) u_0(\tilde{B}_t^x) \exp\left(4\lambda^2 X^m\right) \right]$$

$$= \mathbb{E} \left[ u_0(B_t^x) u_0(\tilde{B}_t^x) \exp\left(4\lambda^2 \left(A_{t,B^x}^{\varepsilon, \delta}, A_{t,\tilde{B}^x}^{\varepsilon, \delta}\right)_{\mathcal{H}}\right) \right].$$

From the boundedness of $u_0$ and [17, Eq. (3.27)] we eventually obtain

$$\sup_{0 < \varepsilon, \delta \le 1} \left\| u_{t,x}^{\varepsilon, \delta} \right\|_{\lambda} < \infty,$$

which shows $u_{t,x} \in \mathcal{G}_{\lambda}$. The estimate (34) follows from the weak lower semi-continuity of the norm. □

### 5.2.2 Stratonovich case

In this part we consider the product of $\dot{W}_{t,x}$ and $u_{t,x}$ in (29) in the sense of the Stratonovich integral. We proceed similar as in Theorem 5.10 and use the results made in [17, Sect. 4]. In [17, Sect. 4] a mild solution $(u_{t,x})_{t>0, x\in\mathbb{R}^d}$ of the Stratonovich version of (29), (30) is constructed. As in the Skorokhod case we show that for all $t > 0$, $x \in \mathbb{R}^d$ the random variable $u_{t,x}$ is contained in $\mathcal{G}$. We first state additional assumption on the covariance $\gamma$ and $\Lambda$, see [17, Hypothesis 4.1], and recall some of the results achieved in [17, Sect. 4].

**Assumption 5.11** Let $\gamma$ and $\Lambda$ be given as in Assumption 5.8. Assume additionally that there exists a constant $0 < \beta < 1$ s.t. for any $t \in \mathbb{R}$,

$$0 \le \gamma(t) \le C_\beta |t|^{-\beta}$$

for some constant $0 < C_\beta < \infty$ and the measure $\sigma$ satisfies

$$\int_{\mathbb{R}^d} \frac{1}{1 + |\xi|^{2-2\beta}} \sigma(d\xi) < \infty.$$

The candidate solution $(u_{t,x})_{t>0, x\in\mathbb{R}^d}$ of (29), (30) for the Stratonovich case is given as a limit of approximations $(u_{t,x}^{\varepsilon, \delta})_{t>0, x\in\mathbb{R}^d}$ where $\varepsilon, \delta > 0$ are cut off parameters and tend to zero. The convergence takes place in $L^p(\mathcal{N}', \mu)$ for every $p \in [1, \infty)$ and uniformly in $t > 0$, $x \in \mathbb{R}^d$, see [17, Proposition 4.7]. The approximations $(u_{t,x}^{\varepsilon,\delta})_{t>0,x\in\mathbb{R}^d}$ are given as a Bochner integral in $L^p(\mathcal{N}', \mu)$, $p \in [1, \infty)$ as above, see [17, Equation (4.15)]. Indeed, for $t > 0$ and $x \in \mathbb{R}^d$ $u_{t,x}^{\varepsilon,\delta}$ is given as

$$
\begin{aligned}
u_{t,x}^{\varepsilon,\delta} &= \mathbb{E}\left[ u_0(B_t^x) \exp\left( A_{t,B^x}^{\varepsilon,\delta} \right) \right] \\
&= \mathbb{E}\left[ u_0(B_t^x) \exp\left( \frac{1}{2} \left\| A_{t,B^x}^{\varepsilon,\delta} \right\|_{\mathcal{H}}^2 \right) : \exp\left( \left\langle A_{t,B^x}^{\varepsilon,\delta}, \cdot \right\rangle \right) : \right].
\end{aligned}
$$

Now we can formulate our result regarding the process $(u_{t,x})_{t>0,x\in\mathbb{R}^d}$.

**Theorem 5.12** *Let $u_0 : \mathbb{R}^d \longrightarrow \mathbb{R}$ be continuous and bounded and $\gamma$ and $\Gamma$ satisfy Assumption 5.11. For every $t > 0$, $x \in \mathbb{R}^d$ the element $u_{t,x} = \lim_{\varepsilon\to 0} \lim_{\delta\to 0} u_{t,x}^{\varepsilon,\delta}$ is contained in $\mathcal{G}$ and for $\lambda \in (0, \infty)$ the $\mathcal{G}_\lambda$-norm of $u_{t,x}$ can be estimated by*

$$
\begin{aligned}
\| u_{t,x} \|_\lambda^2 &\leq \mathbb{E}\Bigg[ u_0(B_t^x) u_0(\tilde{B}_t^x) \exp\Bigg( \frac{1}{2} \Bigg( \int_0^t \int_0^t \gamma(r-s) \Lambda(B_r^x - B_s^x) \, ds \, dr \\
&\quad + \int_0^t \int_0^t \gamma(r-s) \Lambda(\tilde{B}_r^x - \tilde{B}_s^x) \, ds \, dr \Bigg) \Bigg) \\
&\quad \times \exp\Bigg( 4\lambda^2 \int_0^t \int_0^t \gamma(r-s) \Lambda(B_r^x - \tilde{B}_s^x) \, ds \, dr \Bigg) \Bigg],
\end{aligned}
$$

*where $\mathbb{E}$ denotes integration w.r.t. $\mathbb{P}$.*

**Proof** We show this results as in the proof of Theorem 5.10. Indeed, the same arguments as above lead to the fact that for every $\varepsilon, \delta, \lambda, t > 0$ and $x \in \mathbb{R}^d$ the $\mathcal{G}_\lambda$-norm of $u_{t,x}^{\varepsilon,\delta}$ is given by

$$
\begin{aligned}
\left\| u_{t,x}^{\varepsilon,\delta} \right\|_\lambda^2 &= \mathbb{E}\Bigg[ u_0(B_t^x) u_0(\widetilde{B}_t^x) \exp\Bigg( \frac{1}{2} \Bigg( \left\| A_{t,B^x}^{\varepsilon,\delta} \right\|_{\mathcal{H}}^2 + \left\| A_{t,\tilde{B}^x}^{\varepsilon,\delta} \right\|_{\mathcal{H}}^2 \Bigg) \\
&\quad + 4\lambda^2 \left( A_{t,B^x}^{\varepsilon,\delta}, A_{t,\tilde{B}^x}^{\varepsilon,\delta} \right)_{\mathcal{H}} \Bigg) \Bigg].
\end{aligned}
$$

Using [17, Eq. (4.17) ff.] we conclude that

$$
\sup_{0<\varepsilon,\delta\leq 1} \left\| u_{t,x}^{\varepsilon,\delta} \right\|_\lambda < \infty,
$$

implying that $u_{t,x} \in \mathcal{G}_\lambda$ for all $\lambda \in (0, \infty)$. The last part of the statement follows by the same argument as in 5.10 and [17, Eq. (4.6)]. □

**Remark 5.13** We want to point out that the results of Theorems 5.10 and 5.12 give some additional insight into the solution of the stochastic heat equation (29), (30) for the Skorokhod and Stratonovich case. In particular, in both cases the random variable $u_{t,x}$, $t > 0$, $x \in \mathbb{R}^d$, is contained in $\mathcal{G}$. This implies that $u_{t,x}$ is infinitely often Malliavin differentiable and the derivatives of arbitrary order are integrable of order $p$, where $p \in [1, \infty)$ can be arbitrarily large, see [25]. As far as the authors know, this has not been shown for this general class of covariances.

## 6 Outlook: application to stochastic currents

The concept of current is fundamental in geometric measure theory. The simplest version of current is given by the functional

$$\varphi \mapsto \int_0^T \big( \varphi(\gamma(t)), \gamma'(t) \big)_{\mathbb{R}^d} \, dt, \quad 0 < T < \infty,$$

where $\varphi : \mathbb{R}^d \to \mathbb{R}^d$, $d \in \mathbb{N}$, and $\gamma : [0, T] \to \mathbb{R}^d$ is a rectifiable curve. Its vector valued integral kernel informally is given by

$$\zeta(x) = \int_0^T \delta(x - \gamma(t)) \gamma'(t) \, dt, \quad x \in \mathbb{R}^d,$$

where $\delta$ is the Dirac delta. The interested reader may find comprehensive account on the subject in the books [7,21].

A stochastic analog of the current $\zeta$ arises if we replace the deterministic curve $\gamma$ for example by the trajectory of a Brownian motion $(B(t))_{0 \le t \le T}$ taking values in $\mathbb{R}^d$. In this way, we obtain the following functional

$$\xi(x) := \int_0^T \delta(x - B(t)) \, dB(t), \quad x \in \mathbb{R}^d. \tag{35}$$

In the forthcoming manuscript [5] a rigorous definition of (35) is given. Using Wick products the stochastic integral w.r.t. $\mathbb{R}^d$-valued Brownian motions can be defined in the space of Hida distributions. Then our improved characterization of $\mathcal{G}'$ is applied to analyze the regularity of $\xi(x)$, $x \in \mathbb{R}^d$.

There have been some other approaches to study stochastic current, such as Malliavin calculus and stochastic integrals via regularization, see [8–10,13], among others. In [9] $\xi$ was constructed in a negative Sobolev space, i.e., in a generalized function space in the variable $x \in \mathbb{R}^d$. Then the constructed distribution was applied to a model of random vortex filaments in turbulent fluids. The construction in [5] gives for the same object a rigorous definition pointwise in $x \in \mathbb{R}^d \backslash \{0\}$.

**Funding** Open Access funding enabled and organized by Projekt DEAL.

**Publisher's Note** Springer Nature remains neutral with regard to jurisdictional claims in published maps and institutional affiliations.